\documentclass[11pt]{amsart}

\usepackage{amsmath, amssymb, amsthm, amsfonts, amsxtra, mathrsfs}

\input xy
\xyoption{all}

\usepackage{hyperref}

\swapnumbers
\numberwithin{equation}{section}

\theoremstyle{plain}
\newtheorem{theorem}[subsubsection]{Theorem}

\theoremstyle{definition}

\newtheorem{exam}[subsubsection]{Example}

\def\AA{\mathbb{A}}

\def\CC{\mathbb{C}}

\def\FF{\mathbb{F}}
\def\GG{\mathbb{G}}

\def\OO{\mathbb{O}}

\def\QQ{\mathbb{Q}}

\def\ZZ{\mathbb{Z}}

\newcommand\cA{\mathcal{A}}
\newcommand\cB{\mathcal{B}}

\newcommand\cE{\mathcal{E}}
\newcommand\cF{\mathcal{F}}

\newcommand\cH{\mathcal{H}}

\newcommand\cL{\mathcal{L}}
\newcommand\cM{\mathcal{M}}
\newcommand\cN{\mathcal{N}}
\newcommand\cO{\mathcal{O}}
\newcommand\cP{\mathcal{P}}

\newcommand\cZ{\mathcal{Z}}

\newcommand{\sH}{\mathscr{H}}

\newcommand{\sL}{\mathscr{L}}

\def\bR{\mathbf{R}}

\newcommand\frC{\mathfrak{C}}

\newcommand\frX{\mathfrak{X}}

\newcommand\frp{\mathfrak{p}}
\newcommand\frq{\mathfrak{q}}

\newcommand\fru{\mathfrak{u}}

\newcommand\dG{\widehat{G}}
\newcommand\dH{\widehat{H}}

\newcommand{\Bun}{\textup{Bun}}

\newcommand{\coker}{\textup{coker}}

\newcommand\diag{\textup{diag}}

\renewcommand\div{\textup{div}}

\newcommand\Fr{\textup{Fr}}

\newcommand\Frob{\textup{Frob}}
\newcommand\Gal{\textup{Gal}}

\newcommand{\Gr}{\textup{Gr}}

\newcommand{\Hk}{\textup{Hk}}

\newcommand\IC{\textup{IC}}
\newcommand\id{\textup{id}}

\newcommand\inv{\textup{inv}}

\newcommand\loc{\textup{loc}}
\newcommand\Mat{\textup{Mat}}

\newcommand{\Nm}{\textup{Nm}}

\newcommand{\Pic}{\textup{Pic}}

\newcommand{\Res}{\textup{Res}}

\newcommand\RTr{\textup{RTr}}
\newcommand\Sht{\textup{Sht}}

\newcommand\Spec{\textup{Spec}\ }
\newcommand\Spf{\textup{Spf}\ }

\newcommand{\Tr}{\textup{Tr}}

\newcommand{\vol}{\textup{vol}}

\newcommand\Hom{\textup{Hom}}

\newcommand\GL{\textup{GL}}
\newcommand\gl{\mathfrak{gl}}
\newcommand\PGL{\textup{PGL}}
\newcommand\Sg{\textup{S}}
\newcommand\Ug{\textup{U}}

\newcommand{\Gm}{\GG_m}

\newcommand{\Ad}{\textup{Ad}}

\newcommand{\incl}{\hookrightarrow}
\newcommand{\isom}{\stackrel{\sim}{\to}}
\newcommand{\dto}{\dashrightarrow}

\newcommand{\Ql}{\QQ_{\ell}}
\newcommand{\Qlbar}{\overline{\QQ}_\ell}

\newcommand{\hotimes}{\widehat{\otimes}}

\newcommand{\jiao}[1]{\langle{#1}\rangle}
\newcommand{\wt}[1]{\widetilde{#1}}

\newcommand\quash[1]{}
\newcommand\mat[4]{\left(\begin{array}{cc} #1 & #2 \\ #3 & #4 \end{array}\right)}  
\newcommand\un{\underline}

\newcommand{\ov}{\overline}
\newcommand{\bs}{\backslash}

\newcommand\sss{\subsubsection}
\newcommand\xr{\xrightarrow}
\newcommand\op{\oplus}
\newcommand\ot{\otimes}
\newcommand\one{\mathbf{1}}

\newcommand{\cohog}[2]{\textup{H}^{#1}({#2})}     
\newcommand{\cohoc}[2]{\textup{H}_{c}^{#1}({#2})}     
\newcommand{\IHc}[2]{\textup{IH}_{c}^{#1}({#2})}  

\renewcommand\a\alpha
\renewcommand\b\beta
\newcommand\g\gamma
\renewcommand\d\delta
\newcommand\D\Delta

\renewcommand{\th}{\theta}

\newcommand{\ph}{\varphi}
\newcommand{\io}{\iota}
\renewcommand{\r}{\rho}
\renewcommand{\k}{\kappa}
\newcommand{\s}{\sigma}
\newcommand{\Sig}{\Sigma}

\newcommand{\y}{\eta}

\newcommand{\vp}{\varpi}
\renewcommand{\l}{\lambda}
\renewcommand{\L}{\Lambda}
\newcommand{\om}{\omega}

\newcommand\ds{\diamondsuit}

\newcommand{\kbar}{\overline{k}}

\newcommand{\Si}{\Sigma_{\infty}}
\newcommand{\Eis}{\textup{Eis}}

\title{Hitchin type moduli stacks in automorphic representation theory}

\author{Zhiwei Yun}

\address{Department of Mathematics, MIT, 77 Massachusetts Ave, Cambridge, MA 02139}
\email{zyun@mit.edu}

\begin{document}

\begin{abstract}
This is an article prepared for the proceedings of the ICM 2018. In the study of automorphic representations over a function field, Hitchin moduli stack and its variants naturally appear and their geometry helps the comparison of trace formulae.   We give a survey on applications of this observation to a relative fundamental lemma, the arithmetic fundamental lemma and to the higher Gross-Zagier formula.
\end{abstract}

\maketitle

\section{Introduction}

\subsection{Hitchin's original construction}\label{ss:intro Hit}
In an influential paper by Hitchin \cite{Hit}, he introduced the famous integrable system, the moduli space of Higgs bundles. Let $X$ be a smooth proper and geometrically connected curve over a field $k$. Let $G$ be a connected reductive group over $k$. Let $\cL$ be a line bundle over $X$. An {\em $\cL$-twisted $G$-Higgs bundle} over $X$ is a pair $(\cE,\ph)$ where $\cE$ is a principal $G$-bundle over $X$ and $\ph$ is a global section of the vector bundle $\Ad(\cE)\otimes\cL$ over $X$. Here, $\Ad(\cE)$ is the vector bundle associated to $\cE$ and the adjoint representation of $G$. The moduli stack $\cM_{G,\cL}$ of $\cL$-twisted Higgs $G$-bundles over $X$ is called the Hitchin moduli stack. Hitchin defined a map 
\begin{equation}\label{intro Hit fib}
f: \cM_{G,\cL}\to \cA_{G,\cL}
\end{equation}
to some affine space $\cA_{G,\cL}$ by collecting invariants of $\ph$ such as its trace and determinant in the case $G=\GL_{n}$. The map $f$ is called the {\em Hitchin fibration}. When $\cL=\om_{X}$ is the line bundle of $1$-forms on $X$, Hitchin showed that $f$ exhibited the stable part of $\cM_{G,\om_{X}}$ as a completely integrable system. He also gave concrete descriptions of the fibers of $f$ in terms of spectral curves.

\subsection{Applications in geometric representation theory} Although discovered in the context of Yang-Mills theory, Hitchin moduli stacks have subsequently played important roles in the development of geometric representation theory. 

When $\cL=\om_{X}$, $\cM_{G,\om_{X}}$ is essentially the total space of the cotangent bundle of the moduli stack $\Bun_{G}$ of $G$-bundles over $X$. Therefore the categories of twisted $D$-modules on $\Bun_{G}$ give quantizations of $\cM_{G,\om_{X}}$. Beilinson and Drinfeld studied such quantizations and used them to realize part of the geometric Langlands correspondence (namely when the $\dG$-connection comes from an oper). This can be viewed as a global analogue of the Beilinson-Bernstein localization theorem. A related construction in positive characteristic was initiated by Bezrukavnikov and Braverman \cite{BB} for $\GL_{n}$ and extended to any reductive $G$ by T-H.Chen and X.Zhu \cite{CZ}. Hitchin's integrable system also plays a key role in the work of Kapustin and Witten \cite{KW} which ties geometric Langlands correspondence to quantum field theory.

Hitchin moduli stacks have also been used to construct representations of the double affine Hecke algebra, giving global analogues of Springer representations. See \cite{YGS}, \cite{YGSLD} and \cite{OY}.

\subsection{Applications in automorphic representation theory}
B.C. Ng\^o \cite{NgoHit} made the fundamental observation that point-counting on Hitchin fibers is closely related to orbital integrals that appear in the study of automorphic representations.   This observation, along with ingenious technical work, allowed Ng\^o to prove the Lie algebra version of the Fundamental Lemma conjectured by Langlands and Shelstad in the function field case, see \cite{NgoFL}. 

Group versions of the Hitchin moduli stack were introduced by Ng\^o using Vinberg semigroups. They are directly related to the Arthur-Selberg trace formula, as we will briefly review in Section \ref{ss:TF}. See recent works by \cite{NL}, \cite{Bouthier} and \cite{BNS} for applications of group versions of the Hitchin moduli stack.

\subsection{Contents of this report} This report will focus on further applications of variants of Hitchin moduli stacks to automorphic representation theory.

In Section \ref{s:TF}, we explain, in heuristic terms, why Hitchin-type moduli stacks naturally show up in the study of Arthur-Selberg trace formula and more generally, relative trace formulae. A relative trace formula calculates the $L^{2}$-pairing of two distributions on the space of automorphic forms of $G$ given by two subgroups. Such pairings, when restricted to cuspidal automorphic representations, are often related to {\em special values} of automorphic $L$-functions. In Section \ref{ss:JR}, we will elaborate on the relative trace formulae introduced by Jacquet and Rallis, for which the fundamental lemma was proved in \cite{YJR}.

In Section \ref{s:Sht}, we point out a new direction initiated in the works \cite{YAFL}, \cite{YZ1} and \cite{YZ2}. Drinfeld introduced the moduli stack of Shtukas as an analogue of Shimura varieties for function fields, which turns out to allow richer variants than Shimura varieties. Cohomology classes of these moduli of Shtukas generalize the notion of automorphic forms. In Section \ref{ss:Sht} we review the basic definitions of Shtukas, and discuss the spectral decomposition for the cohomology of moduli of Shtukas. In Section \ref{ss:HD}, we introduce Heegner-Drinfeld cycles on the moduli of $G$-Shtukas coming from subgroups $H$ of $G$. The relative trace in the context of Shtukas is then defined in Section \ref{ss:Sht RTF} as the intersection pairing of two Heegner-Drinfeld cycles. We believe that such pairings, when restricted to the isotypical component of a cuspidal automorphic representation, are often related to {\em higher derivatives} of automorphic $L$-functions.  We then explain why Hitchin-type moduli stacks continue to play a key role in the Shtuka context, and what new geometric ingredients are needed to study relative trace formulae in this setting. 

In Section \ref{ss:YZ}-\ref{ss:AFL} we survey what has been proven in this new direction. In Section \ref{ss:YZ}, we review \cite{YZ1} and \cite{YZ2}, in which we obtain formulae relating higher derivatives of automorphic $L$-functions for $\PGL_{2}$ and the intersection numbers of Heegner-Drinfeld cycles. Our results generalize the Gross-Zagier formula in the function field case.  In Section \ref{ss:AFL} we discuss the analogue of the fundamental lemma in the Shtuka setting. This was originally conjectured by W.Zhang under the name {\em arithmetic fundamental lemma}. We state an extension of his conjecture for function fields involving higher derivatives of orbital integrals, and sketch our strategy to prove it.

\subsection*{Acknowledgments} 
I would like to thank my mentors for bringing me to the exciting crossroads of algebraic geometry, representation theory and number theory. The influences of R. MacPherson, P. Deligne, G. Lusztig, B. Gross, B.C. Ng\^o and R. Bezrukavnikov on my career are especially important. I also learned a lot from my collaborators, colleagues and fellow students, to whom I would like to express my deep gratitude.

\section{Hitchin moduli stack and trace formulae}\label{s:TF}

Throughout this article we fix a finite field $k=\FF_{q}$. Let $X$ be a smooth, projective and geometrically connected curve over $k$ of genus $g$. Let $F=k(X)$ be the function field of $X$. The places of $F$ can be identified with the set $|X|$ of closed points of $X$. Let $\AA$ denote the ring of ad\`eles of $F$. For $x\in |X|$, let $\cO_{x}$ be the completed local ring of $X$ at $x$, and $F_{x}$ (resp. $k_{x}$) be its fraction field (resp. residue field). We also use $\vp_{x}$ to denote a uniformizer of $\cO_{x}$.

In this section we work with the classical notion of automorphic forms for groups over function fields. We shall briefly review the Arthur-Selberg trace formula and the relative trace formulae, and explain why Hitchin type moduli stacks naturally show up in the study of these trace formulae.

\subsection{Arthur-Selberg trace formula}\label{ss:TF} 
The Arthur-Selberg trace formula is an important tool in the theory of automorphic representations. For a detailed introduction to the theory over a number field we recommend Arthur's article \cite{Ar}.  Here we focus on the function field case. The idea that Hitchin moduli stacks give a geometric interpretation of Arthur-Selberg trace formula is due to B.C. Ng\^o. For more details, we refer to \cite{NgoHit} for the Lie algebra version, and \cite{FN} for the group version.

We ignore the issue of convergence in the discussion (i.e., we pretend to be working with an anisotropic group $G$), but we remark that the convergence issue lies at the heart of the theory of Arthur and we are just interpreting the easy part of his theory from a geometric perspective.

\sss{The classical setup} Let $G$ be a split connected reductive group over $k$ and we view it as a group scheme over $X$ (hence over $F$) by base change. Automorphic forms for $G$ are $\CC$-valued smooth functions on the coset space $G(F)\bs G(\AA)$. Fix a Haar measure $\mu_{G}$ on $G(\AA)$. Let $\cA$ be the space of automorphic forms for $G$. For any smooth compactly supported function $f$ on $G(\AA)$, it acts on $\cA$ by right convolution $R(f)$.

The Arthur-Selberg trace formula aims to express the trace of $R(f)$ on $\cA$ in two different ways: one as a sum over conjugacy classes of $G(F)$ (the geometric expansion) and the other as a sum over automorphic representations (the spectral expansion). The primitive form of the geometric expansion reads 
\begin{equation}\label{geom exp TF}
\Tr(R(f), \cA)''=''\sum_{\g\in G(F)/\sim} J_{\g}(f)
\end{equation}
where $\g$ runs over $G(F)$-conjugacy classes in $G(F)$, and $J_{\g}(f)$ is the orbital integral
\begin{equation*}
J_{\g}(f)=\vol(G_{\g}(F)\bs G_{\g}(\AA), \mu_{G_{\g}})\int_{G_{\g}(\AA)\bs G(\AA)}f(g^{-1}\g g)\frac{\mu_{G}}{\mu_{G_{\g}}}(g)
\end{equation*}
where $\mu_{G_{\g}}$ is any Haar measure on the centralizer $G_{\g}(\AA)$ of $\g$. We write the equality sign in quotation marks \footnote{Later we will use the same notation (equal signs in quotation marks) to indicate  heuristic identities.} to indicate that the convergence issue has been ignored in \eqref{geom exp TF}. We will give a geometric interpretation of the geometric expansion.

Fix a compact open subgroup $K=\prod_{x\in |X|}K_{x}\subset G(\AA)$ and assume $\vol(K,\mu_{G})=1$. Let $\cA_{K}=C_{c}(G(F)\bs G(\AA)/K)$ on which the Hecke algebra $C_{c}(K\bs G(\AA)/K)$ acts. For $g\in G(\AA)$, there is a Hecke correspondence attached to the double coset $KgK\subset G(\AA)$
\begin{equation}\label{Hk classical}
\xymatrix{G(F)\bs G(\AA)/K & G(F)\bs G(\AA)/(K\cap gKg^{-1})\ar[l]_-{p_{0}}\ar[r]^-{q_{0}} & G(F)\bs G(\AA)/K
}
\end{equation}
where $p_{0}$ is the natural projection and $q_{0}$ is induced by right multiplication by $g$. The action of $f=\one_{KgK}$ on $\cA_{K}$ is given by $\ph\mapsto q_{0!}p^{*}_{0}\ph$, where $q_{0!}$ means summing over the fibers of $q_{0}$. Upon ignoring convergence issues, the trace of $R(\one_{KgK})$ on $\cA_{K}$ is equal to the cardinality of the restriction of $G(F)\bs G(\AA)/(K\cap gKg^{-1})$ to the diagonal $G(F)\bs G(\AA)/K$ via the maps $(p_{0},q_{0})$. In other words, we should form the pullback diagram of groupoids
\begin{equation}\label{MHk cl}
\xymatrix{M_{G,KgK} \ar[d]\ar[r] & G(F)\bs G(\AA)/(K\cap gKg^{-1})\ar[d]^{(p_{0},q_{0})}\\
G(F)\bs G(\AA)/K\ar[r]^-{\D} & G(F)\bs G(\AA)/K\times G(F)\bs G(\AA)/K}
\end{equation}
and we have a heuristic identity
\begin{equation}\label{TF cl}
\Tr(R(\one_{KgK}), \cA_{K})''=''\#M_{G,KgK}.
\end{equation}
Here $\#\frX$ of a groupoid $\frX$ is a counting of isomorphism classes of objects in $\frX$ weighted by the reciprocals of the sizes of automorphism groups.

\sss{Weil's interpretation} Let $K_{0}=\prod_{x\in |X|}G(\cO_{x})$.  It was observed by Weil that the double coset groupoid $G(F)\bs G(\AA)/K_{0}$ is naturally isomorphic to the groupoid of $G$-bundles over $X$ that are trivial at the generic point of $X$. In fact, starting from $g=(g_{x})\in G(\AA)$, one assigns the $G$-bundle on $X$ that is glued from the trivial bundles on the generic point $\Spec F$ and the formal disks $\Spec \cO_{x}$ via the ``transition matrices'' $g_{x}$. For a compact open $K\subset K_{0}$,  one can similarly interpret $G(F)\bs G(\AA)/K$ as the groupoid of $G$-bundles with $K$-level structures. There is an algebraic stack $\Bun_{G,K}$ classifying  $G$-bundles on $X$ with $K$-level structures, and the above observation can be rephrased as a fully faithful embedding of groupoids
\begin{equation}\label{equiv Bun}
G(F)\bs G(\AA)/K\incl \Bun_{G,K}(k).
\end{equation}
A priori, the groupoid $\Bun_{G,K}(k)$ contains also $G$-bundles that are not trivial at the generic point, or equivalently $G'$-bundles for certain inner forms $G'$ of $G$. Since we assume $G$ is split, the embedding \eqref{equiv Bun} is in fact an equivalence.

In the same spirit, we interpret $G(F)\bs G(\AA)/(K\cap gKg^{-1})$ as the groupoid of triples $(\cE, \cE', \a)$ where $\cE,\cE'$ are $G$-bundles with $K$-level structures on $X$, and $\a: \cE\dto\cE'$ is a {\em rational} isomorphism between $\cE$ and $\cE'$ (i.e., an isomorphism of $G$-bundles over the generic point $\Spec F$) such that the relative position of $\cE$ and $\cE'$ at each closed point $x\in |X|$ is given by the double coset $K_{x}g_{x}K_{x}$. For example, when $G=\GL_{n}$, $K_{x}=\GL_{n}(\cO_{x})$ and $g_{x}=\diag(\vp_{x}, 1,\cdots, 1)$, then $\a$ has relative position $K_{x}g_{x}K_{x}$ at $x$ if and only if $\a$ extends to a homomorphism $\a_{x}: \cE|_{\Spec\cO_{x}}\to \cE'|_{\Spec \cO_{x}}$, and that $\coker(\a_{x})$ is one-dimensional over the residue field $k_{x}$. There is a moduli stack $\Hk_{G,KgK}$ classifying such triples $(\cE, \cE', \a)$. The above discussion can be  rephrased as an equivalence of groupoids
\begin{equation}\label{equiv Hk}
\Hk_{G,KgK}(k)\cong G(F)\bs G(\AA)/(K\cap gKg^{-1}).
\end{equation}
Moreover, $\Hk_{G,KgK}$ is equipped with two maps to $\Bun_{G,K}$ by recording $\cE$ and $\cE'$, which allow us to view it as a self-correspondence of $\Bun_{G,K}$
\begin{equation}\label{Hk geometric}
\xymatrix{\Bun_{G,K} & \Hk_{G,KgK}\ar[l]_-{p}\ar[r]^-{q}  & \Bun_{G,K}}
\end{equation}
Under the equivalences \eqref{equiv Bun} and \eqref{equiv Hk}, the diagram \eqref{Hk geometric} becomes the diagram \eqref{Hk classical} after taking $k$-points.

\sss{Geometric interpretation of the trace} Continuing further with Weil's observation, we can form the stack-theoretic version of \eqref{MHk cl}, and define the stack $\cM_{G,KgK}$ by the Cartesian diagram 
\begin{equation}\label{MHk geom}
\xymatrix{\cM_{G,KgK} \ar[d]\ar[r] & \Hk_{G,KgK}\ar[d]^{(p,q)}\\
\Bun_{G,K}\ar[r]^-{\D} & \Bun_{G,K}\times \Bun_{G,K}}
\end{equation}
so that  $\cM_{G,KgK}(k)=M_{G,KgK}$. 

By the defining Cartesian diagram \eqref{MHk geom}, $\cM_{G,KgK}$ classifies pairs $(\cE,\ph)$ where $\cE$ is a $G$-bundle over $X$ with $K$-level structures, and $\ph: \cE\dto \cE$ is a rational {\em automorphism} with relative position given by $KgK$. 

Recall the classical Hitchin moduli stack $\cM_{G,\cL}$ in  Section \ref{ss:intro Hit}. If we write $\cL=\cO_{X}(D)$ for some effective divisor $D$, $\cM_{G,\cL}$ then classifies pairs $(\cE, \ph)$ where $\cE$ is a $G$-bundle over $X$ and $\ph$ is a rational section of $\Ad(\cE)$ (an infinitesimal automorphism of $\cE$) with poles controlled by $D$. Therefore $\cM_{G,\cL}$ can be viewed as a Lie algebra version of $\cM_{G,KgK}$, and $\cM_{G,KgK}$ is a group version of $\cM_{G,\cL}$.

Let $\frC_{G}$ be the GIT quotient of $G$ by the conjugation action of $G$;  $\frC_{G}(F)$ is the set of {\em stable conjugacy classes} in $G(F)$. 
There is an affine scheme $\cB_{G,KgK}$ classifying rational maps $X\dto \frC_{G}$ with poles controlled by $KgK$. The Hitchin fibration \eqref{intro Hit fib} has an analogue
\begin{equation}\label{Hitchin map}
h_{G}:\cM_{G,KgK}\to \cB_{G,KgK}.
\end{equation}
Using the map $h_{G}$, the counting of $M_{G,KgK}$, hence the trace of $R(\one_{KgK})$, can be decomposed into a sum over certain stable conjugacy classes $a$
\begin{equation}\label{TF over a}
\Tr(R(\one_{KgK}), \cA_{K})''=''\sum_{a\in \cB_{G,KgK}(k)}\# \cM_{G,KgK}(a)(k).
\end{equation}
Here $\cM_{G,KgK}(a)$ (a stack over $k$) is the fiber $h_{G}^{-1}(a)$.  To tie back to the classical story,  $\#\cM_{G,KgK}(a)(k)$ is in fact a sum of orbital integrals
\begin{equation*}
\# \cM_{G,KgK}(a)(k)=\sum_{\g\in G(F)/\sim, [\g]=a} J_{\g}(\one_{KgK})
\end{equation*}
over $G(F)$-conjugacy classes $\g$ that belong to the stable conjugacy class $a$. 

By the Lefschetz trace formula, we can rewrite \eqref{TF over a} as
\begin{equation}\label{TF coho}
\Tr(R(\one_{KgK}), \cA_{K})''=''\sum_{a\in \cB_{G,KgK}(k)}\Tr(\Frob_{a}, (\bR h_{G!}\Ql)_{a}).
\end{equation}
This formula relates the Arthur-Selberg trace to the direct image complex of the Hitchin fibration $h_{G}$ (called the {\em Hitchin complex} for $G$). Although it is still difficult to get a closed formula for each term in \eqref{TF coho}, this geometric point of view can be powerful in comparing traces for two different groups $G$ and $H$ by relating their Hitchin bases and Hitchin complexes. 

In the work of Ng\^o \cite{NgoFL}, where the Lie algebra version was considered, the Hitchin complex was studied in depth using tools such as perverse sheaves and the decomposition theorem. When $H$ is an endoscopic group of $G$, Ng\^o shows that the stable part of the Hitchin complex for $H$ appears as a direct summand of the Hitchin complex for $G$, from which he deduces the Langlands-Shelstad fundamental lemma for Lie algebras over function fields.

\subsection{Relative trace formulae}\label{ss:RTF}

\sss{Periods of automorphic forms}
For simplicity we assume $G$ is semisimple. Let $H\subset G$ be a subgroup defined over $F$, and $\mu_{H}$ a Haar measure on $H(\AA)$. For a cuspidal automorphic representation $\pi$ of $G(\AA)$, the linear functional on $\pi$
\begin{eqnarray}\label{per}
\cP^{G}_{H,\pi}&:& \pi\to \CC\\
\notag&& \ph\mapsto \int_{H(F)\bs H(\AA)}\ph\mu_{H}
\end{eqnarray}
is called the {\em $H$-period} of  $\pi$. It factors through the space of coinvariants $\pi_{H(\AA)}$. 

One can also consider variants  where we integrate $\ph$ against an automorphic character $\chi$ of $H(\AA)$. If $\pi$ has nonzero $H$-period, it is called $H$-distinguished. Distinguished representations are used to characterize important classes of automorphic representations such as those coming by functoriality from another group. In case the local coinvariants $(\pi_{x})_{H(F_{x})}$  are one-dimensional for almost all places $x$ (as in the case for many spherical subgroups of $G$), one expects the period $\cP^{G}_{H,\pi}$ to be related to special values of $L$-functions of $\pi$. 

\begin{exam} Let $G=\PGL_{2}$ and $H=A$ be the diagonal torus. Then by Hecke's theory, for a suitably chosen $\ph\in \pi$, $\int_{H(F)\bs H(\AA)}\ph(t)|t|^{s-1/2}dt$ gives the standard $L$-function $L(s, \pi)$. 
\end{exam}

\sss{Relative trace formulae} Now suppose $H_{1}, H_{2}$ are two subgroups of $G$. Let $\pi$ be a cuspidal automorphic representation of $G(\AA)$ and $\wt\pi$ its contragradient. We get a bilinear form 
\begin{equation*}
\cP^{G}_{H_{1},\pi}\ot\cP^{G}_{H_{2},\wt\pi}: \pi\otimes\wt\pi\to  \CC
\end{equation*}
In case the local coinvariants $(\pi_{x})_{H_{i}(F_{x})}$  are one-dimensional for all places $x$ and $i=1,2$, the $H_{1}(\AA)\times H_{2}(\AA)$ invariant bilinear forms on $\pi\times\wt\pi$ are unique up to scalar. Therefore $\cP^{G}_{H_{1},\pi}\ot\cP^{G}_{H_{2},\wt\pi}$ is a multiple of the natural pairing on $\pi\times \wt\pi$ given by the Petersson inner product. This multiple is often related to special values of $L$-functions attached to $\pi$. For a systematic treatment of this topic, see the book by Sakellaridis and Venkatesh \cite{SV}.

An important tool to study the pairing $\cP^{G}_{H_{1},\pi}\ot\cP^{G}_{H_{2},\wt\pi}$ is the relative trace formula. We have natural  maps of cosets 
\begin{equation}\label{phH auto}
\xymatrix{H_{1}(F)\bs H_{1}(\AA)\ar[r]^{\ph_{1}} & G(F)\bs G(\AA) & H_{2}(F)\bs H_{2}(\AA)\ar[l]_{\ph_{2}}}
\end{equation}
Consider the push-forward of the constant functions on $H_{i}(F)\bs H_{i}(\AA)$ along $\ph_{i}$, viewed as distributions on $G(F)\bs G(\AA)$. Since we will only give a heuristic discussion of the relative trace formula, we will pretend that the $L^{2}$-pairing of two distributions makes sense. The relative trace of a test function $f\in C^{\infty}_{c}(G(\AA))$ is the $L^{2}$-pairing
\begin{equation}\label{rel trace}
\RTr^{G}_{H_{1}, H_{2}}(f)''=''\jiao{\ph_{1,!}\one_{H_{1}(F)\bs H_{1}(\AA)}, R(f)\ph_{2,!}\one_{H_{2}(F)\bs H_{2}(\AA)}}_{L^{2}(G(F)\bs G(\AA),\mu_{G})}
\end{equation}
A variant of this construction is to replace the constant function $1$ on $H_{i}(F)\bs H_{i}(\AA)$ by an automorphic quasi-character $\chi_{i}$ of $H_{i}(\AA)$. 

The relative trace formula is an equality between a spectral expansion of $\RTr^{G}_{H_{1},H_{2}}(f)$ into quantities related to $\cP^{G}_{H_{1},\pi}\ot\cP^{G}_{H_{2},\wt\pi}$ and a geometric expansion into a sum of orbital integrals $J^{G}_{H_{1},H_{2}, \g}(f)$ indexed by double cosets $H_{1}(F)\g H_{2}(F)\subset G(F)$
\begin{equation}\label{rel orb int}
J^{G}_{H_{1},H_{2}, \g}(f)=\vol(H_{\g}(F)\bs H_{\g}(\AA), \mu_{H_{\g}})\int_{H_{\g}(\AA)\bs (H_{1}\times H_{2})(\AA)}f(h_{1}^{-1}\g h_{2})\frac{\mu_{H_{1}\times H_{2}}}{\mu_{H_{\g}}}(h_{1}, h_{2})
\end{equation}
where $H_{\g}$ is the stabilizer of $\g$ under the left-right translation on $G$ by $H_{1}\times H_{2}$.

Let $\frC_{H_{1}, H_{2}}^{G}=\Spec F[G]^{H_{1}\times H_{2}}$. We may call elements in $\frC_{H_{1}, H_{2}}^{G}(F)$ {\em stable orbits} of $G$ under the action of $H_{1}\times H_{2}$.  There is a tautological map
\begin{equation*}
\inv: H_{1}(F)\bs G(F)/H_{2}(F)\to \frC_{H_{1}, H_{2}}^{G}(F).
\end{equation*}
We define 
\begin{equation}\label{stable orb int}
J^{G}_{H_{1}, H_{2}}(a,f)=\sum_{\g\in H_{1}(F)\bs G(F)/H_{2}(F) , \inv(\g)=a}J^{G}_{H_{1},H_{2}, \g}(f).
\end{equation}

Now we fix a compact open subgroup $K\subset G(\AA)$ and let $K_{i}=K\cap H_{i}(\AA)$ for $i=1,2$. Choose Haar measures on $H_{i}(\AA)$ and $G(\AA)$ so that $K_{i}$ and $K$ have volume $1$. Consider the test function $f=\one_{KgK}$ as before.  Unwinding the definition of the relative trace, we can rewrite \eqref{rel trace} as 
\begin{equation}\label{RTr cl}
\RTr^{G}_{H_{1},H_{2}}(\one_{KgK})''=''\# M^{G}_{H_{1}, H_{2}, KgK}
\end{equation}
where $M^{G}_{H_{1}, H_{2}, KgK}$ is defined by the Cartesian diagram of groupoids
\begin{equation*}
\xymatrix{M^{G}_{H_{1}, H_{2}, KgK}\ar[r]\ar[d] & G(F)\bs G(\AA)/(K\cap gKg^{-1})\ar[d]^{(p_{0},q_{0})}\\
H_{1}(F)\bs H_{1}(\AA)/K_{1}\times H_{2}(F)\bs H_{2}(\AA)/K_{2}\ar[r]^-{(\ph_{1}, \ph_{2})} & G(F)\bs G(\AA)/K\times G(F)\bs G(\AA)/K}
\end{equation*} 

\sss{Geometric interpretation of the relative trace}\label{sss:geom rel trace} Now we give a geometric interpretation of $M^{G}_{H_{1}, H_{2}, KgK}$. We assume $H_{1}$ and $H_{2}$ are also obtained by base change from split reductive groups over $k$, which we denote by the same notation. We have maps 
\begin{equation*}
\xymatrix{\Bun_{H_{1},K_{1}}\ar[r]^{\Phi_{1}} & \Bun_{G,K} & \Bun_{H_{2},K_{2}}\ar[l]_{\Phi
_{2}}.
}
\end{equation*} 
Taking $k$-points of the above diagram we recover \eqref{phH auto} up to modding out by the compact open subgroups $K_{i}$ and $K$. We may form the Cartesian diagram of stacks
\begin{equation*}
\xymatrix{\cM^{G}_{H_{1}, H_{2}, KgK}\ar[r]\ar[d] & \Hk_{G,KgK}\ar[d]^{(p,q)}\\
\Bun_{H_{1},K_{1}}\times \Bun_{H_{2},K_{2}}\ar[r]^-{(\Phi_{1}, \Phi_{2})} & \Bun_{G,K}\times \Bun_{G,K} }
\end{equation*}
We have $\cM^{G}_{H_{1}, H_{2}, KgK}(k)=M^{G}_{H_{1}, H_{2}, KgK}$.

\begin{exam}\label{ex:double} Consider the special case $G=G_{1}\times G_{1}$ ($G_{1}$ is a semisimple group over $k$) and $H_{1}=H_{2}$ is the diagonal copy of $G_{1}$, $K_{1}=K_{2}$, $K=K_{1}\times K_{2}$. Taking $g=(1,g_{1})$ for some $g_{1}\in G_{1}(\AA)$, we get a canonical isomorphism between the stacks $\cM^{G_{1}\times G_{1}}_{\D(G_{1}),\D(G_{1}), K(1,g_{1})K}$ and $\cM_{G_{1},K_{1}g_{1}K_{1}}$. In this case, the relative trace is the usual trace of $R(\one_{K_{1}g_{1}K_{1}})$ on the space of automorphic forms for $G_{1}$.
\end{exam}

The moduli stack $\cM^{G}_{H_{1}, H_{2}, KgK}$ classifies $(\cE_{1}, \cE_{2}, \a)$ where $\cE_{i}$ is an $H_{i}$-bundle with $K_{i}$-structure over $X$ for $i=1,2$; $\a$ is a rational isomorphism between the $G$-bundles induced from $\cE_{1}$ and $\cE_{2}$, with relative position given by $KgK$.

One can construct a scheme $\cB^{G}_{H_{1}, H_{2}, KgK}$ classifying rational maps $X\dto \frC^{G}_{H_{1},H_{2}}$ with poles controlled by $KgK$, so that $\cB^{G}_{H_{1}, H_{2},KgK}(k)\subset \frC^{G}_{H_{1}, H_{2}}(F)$. For  $(\cE_{1}, \cE_{2},\a)\in\cM^{G}_{H_{1}, H_{2}, KgK}$, we may restrict $\a$ to the generic point of $X$ and take its invariants as a rational map $X\dto \frC^{G}_{H_{1},H_{2}}$. This way we get a map of algebraic stacks
\begin{equation}\label{rel Hitchin}
h^{G}_{H_{1}, H_{2}}:\cM^{G}_{H_{1}, H_{2}, KgK}\to \cB^{G}_{H_{1}, H_{2}, KgK}.
\end{equation}
In the situation of Example \ref{ex:double}, $h^{G}_{H_{1}, H_{2}}$ specializes to the usual Hitchin map \eqref{Hitchin map} for $G_{1}$, so we may think of $h^{G}_{H_{1}, H_{2}}$ as an analogue of the Hitchin map for the Hitchin-like moduli $\cM^{G}_{H_{1}, H_{2}, KgK}$. Taking $k$-points of \eqref{rel Hitchin} we get a map
\begin{equation}\label{rel Hitchin cl}
M^{G}_{H_{1}, H_{2}, KgK}\to [H_{1}\bs G/H_{2}](F)\xr{\inv}\frC_{H_{1}, H_{2}}^{G}(F)
\end{equation}
whose fiber over $a\in \frC_{H_{1}, H_{2}}^{G}(F)$ has cardinality equal to $J^{G}_{H_{1},H_{2}}(a, \one_{KgK})$ defined in \eqref{stable orb int}.   We may thus decompose the relative trace into a sum of point-counting along the fibers of the map \eqref{rel Hitchin cl}
\begin{eqnarray}\label{RTr Lef}
\notag \RTr^{G}_{H_{1}, H_{2}}(\one_{KgK})&''=''&\sum_{a\in \frC^{G}_{H_{1},H_{2}}(F)}J^{G}_{H_{1},H_{2}}(a,\one_{KgK})\\
&=&\sum_{a\in \cB^{G}_{H_{1}, H_{2},KgK}(k)} \Tr(\Frob_{a}, (\bR h^{G}_{H_{1}, H_{2},!}\Ql)_{a}).
\end{eqnarray}
The above formula relates the relative trace to the direct image complex $\bR h^{G}_{H_{1}, H_{2},!}\Ql$. As in the case of the Arthur-Selberg trace formula, we may apply sheaf-theoretic tools to study this direct image complex, especially when it comes to comparing two such complexes.

\begin{exam}\label{ex:A}
Consider the case $G=\PGL_{2}$, and $H_{1}=H_{2}=A\subset G$ is the diagonal torus. Let $K_{1}=K_{2}=\prod_{x}A(\cO_{x})$ and $K=\prod_{x}G(\cO_{x})$. Let $D=\sum_{x}n_{x}x$ be an effective divisor. Define the function $h_{D}$ on $G(\AA)$ to be the characteristic function of $\Mat_{2}(\OO)_{D}=\{(g_{x})|g_{x}\in \Mat_{2}(\cO_{x}), v_{x}(\det g_{x})=n_{x}, \forall x\in |X|\}$. Then $\Mat_{2}(\OO)_{D}$ is a finite union of $K$-double cosets $Kg_{i}K$. We will define a stack $\cM^{G}_{A,A,D}$ which turns out to be the union of the $\cM^{G}_{A,A,Kg_{i}K}$.

Consider the stack $\wt\cM^{G}_{A,A,D}$ classifying the data $(\cL_{1}, \cL_{2}, \cL_{1}', \cL_{2}', \ph)$ where
\begin{itemize}
\item $\cL_{i}$, $\cL_{i}'$ are line bundles over $X$, for $i=1,2$;
\item $\ph: \cL_{1}\op \cL_{2}\to \cL_{1}'\op \cL_{2}'$ is an injective map of coherent sheaves such that $\det(\ph)$, viewed as a section of the line bundle $\cL^{-1}_{1}\ot\cL^{-1}_{2}\ot \cL'_{1}\ot \cL_{2}'$, has divisor $D$. 
\end{itemize}
The Picard stack $\Pic_{X}$ acts on $\wt\cM^{G}_{A,A,D}$ by simultaneously tensoring on $\cL_{i}$ and $\cL'_{i}$. We define
\begin{equation*}
\cM^{G}_{A,A,D}\cong \wt\cM^{G}_{A,A,D}/\Pic_{X}.
\end{equation*}
The bi-$A$-invariant regular functions on $G$ are generated by $\mat{a}{b}{c}{d}\mapsto \frac{bc}{ad-bc}$, therefore the space $\frC^{G}_{A,A}$ is isomorphic to $\AA^{1}$. Define the Hitchin base $\cB^{G}_{A,A,D}$ to be the affine space $\cohog{0}{X, \cO_{X}(D)}$.

To define the Hitchin map in this case, we write $\ph$ above as a matrix $\mat{\ph_{11}}{\ph_{12}}{\ph_{21}}{\ph_{22}}$ where $\ph_{ij}$ is a section of $\cL^{-1}_{j}\ot\cL_{i}'$. The determinant $\det(\ph)$ gives an isomorphism $\cL^{-1}_{1}\ot\cL^{-1}_{2}\ot \cL'_{1}\ot \cL_{2}'\cong \cO_{X}(D)$. On the other hand, $\ph_{12}\ph_{21}$ gives another section of $\cL^{-1}_{1}\ot\cL^{-1}_{2}\ot \cL'_{1}\ot \cL_{2}'$. The Hitchin map 
\begin{equation*}
h^{G}_{A, A}: \cM^{G}_{A,A,D}\to \cohog{0}{X, \cO_{X}(D)}=\cB^{G}_{A,A,D}
\end{equation*}
then sends $(\cL_{1}, \cL_{2}, \cL_{1}', \cL_{2}', \ph)$ to $\ph_{12}\ph_{21}$, viewed as a section of $\cO_{X}(D)$ via the identification given by $\det(\ph)$.

Although $\cM^{G}_{A,A,D}$ is not of finite type, it is the disjoint union of finite type substacks indexed by a subset of $\ZZ^{4}/\D(\ZZ)$. Indeed, for  $\un d=(d_{1},d_{2},d_{1}',d_{2}')\in\ZZ^{4}/\D(\ZZ)$ such that $d_{1}'+d_{2}'=d_{1}+d_{2}+\deg D$, the substack ${}^{\un d}\cM^{G}_{A,A,D}$ where $\deg\cL_{i}=d_{i}$ and $\deg \cL_{i}'=d'_{i}$ is of finite type.  We may write $\RTr^{G}_{A,A}(h_{D})$ as a formal sum of 
\begin{equation*}
{}^{\un d}\RTr^{G}_{A,A}(h_{D})=\#{}^{\un d}\cM^{G}_{A,A,D}(k)=\Tr(\Frob, \cohoc{*}{{}^{\un d}\cM^{G}_{A,A,D}\ot\kbar,\Qlbar}).
\end{equation*}
\end{exam}
%

\subsection{Relative fundamental lemma}\label{ss:JR}
In many cases we do not expect to prove closed formulae for relative traces of the form \eqref{rel trace}. Instead, for applications to problems on automorphic representations, it often suffices to establish an identity between relative traces for two different situations $(G,H_{1}, H_{2})$ and $(G',H'_{1}, H'_{2})$.  

\sss{General format of RTF comparison}
In order to establish such an identity, we need the following structures or results:
\begin{enumerate}
\item There should be an isomorphism between the spaces of invariants  $\frC^{G}_{H_{1}, H_{2}}\cong \frC^{G'}_{H'_{1}, H'_{2}}$.

\item (fundamental lemma) For almost all $x\in |X|$, and $a_{x}\in  \frC^{G}_{H_{1}, H_{2}}(F_{x})\cong \frC^{G'}_{H'_{1}, H'_{2}}(F_{x})$, we should have an identity of local orbital integrals up to a transfer factor
\begin{equation*}
J^{G}_{H_{1},H_{2},x}(a_{x},\one_{G(\cO_{x})})\sim  J^{G'}_{H'_{1}, H'_{2},x}(a_{x},\one_{G'(\cO_{x})}).
\end{equation*}
Here $J^{G}_{H_{1},H_{2},x}(a_{x},f_{x})$ is the local analogue of $J^{G}_{H_{1},H_{2}}(a,f)$ defined in \eqref{stable orb int}.

\item (smooth matching) For any $x\in |X|$, $a_{x}\in  \frC^{G}_{H_{1}, H_{2}}(F_{x})\cong \frC^{G'}_{H'_{1}, H'_{2}}(F_{x})$, and $f_{x}\in C_{c}^{\infty}(G(F_{x}))$, there exists $f'_{x}\in C_{c}^{\infty}(G'(F_{x}))$ such that $J^{G}_{H_{1}, H_{2},x}(a_{x},f_{x})=J^{G'}_{H'_{1}, H'_{2}, x}(a_{x},f'_{x})$.
\end{enumerate}
The geometric interpretation \eqref{RTr Lef} of the relative trace gives a way to prove the fundamental lemma by comparing the direct image complexes of the Hitchin maps $h^{G}_{H_{1}, H_{2}}$ and $h^{G'}_{H'_{1}, H'_{2}}$. Below we discuss one such example.

\sss{The relative trace formulae of Jacquet and Rallis}
Jacquet and Rallis \cite{JR} proposed a relative trace formula approach to the Gan-Gross-Prasad conjecture for unitary groups. They considered two relative trace formulae, one involving general linear groups, and the other involving unitary groups. They formulated both the fundamental lemma and the smooth matching in this context as conjectures.
In \cite{YJR}, we used the geometric interpretation sketched in Section \ref{sss:geom rel trace} to prove the fundamental lemma conjectured by Jacquet and Rallis, in the case of function fields. In the appendix to \cite{YJR}, J.Gordon used model theory to deduce the mixed characteristic case from the function field case. On the other hand, W. Zhang \cite{ZFourier} proved the smooth matching for the Jacquet-Rallis relative trace formula at non-archimedean places. Together with the fundamental lemma proved in \cite{YJR}, W.Zhang deduced the Gan-Gross-Prasad conjecture for unitary groups under some local restrictions.

In the next two examples, we introduce the groups involved in the two trace formulae in \cite{JR}, and sketch the definition of the moduli  stacks relevant to the orbital integrals.  Since we proved the fundamental lemma by reducing to its Lie algebra analogue,  our moduli stacks will be linearized versions of the Hitchin-type moduli stacks introduced in Section \ref{sss:geom rel trace}, which are closer to the classical Hitchin moduli stack.

\begin{exam}\label{ex:JR GL} Let $F'/F$ be a separable quadratic extension corresponding to a double cover $\nu: X'\to X$. Let $\s\in\Gal(F'/F)$ be the nontrivial involution. Consider $G=\Res_{F'/F}\GL_{n}\times \Res_{F'/F}\GL_{n-1}$, $H_{1}=\Res_{F'/F}\GL_{n-1}$ and $H_{2}=\GL_{n}\times \GL_{n-1}$ (over $F$). The embedding $H_{1}\to G$ sends $h\in H_{1}$ to $(\mat{h}{}{}{1},h)\in G$.   

The double quotient $H_{1}\bs G/H_{2}$ can be identified with $\GL_{n-1}\bs \Sg_{n}$, where 
$$\Sg_{n}=\{g\in \Res_{F'/F}\GL_{n}|  \s(g)=g^{-1}\}$$  
with $\GL_{n-1}$ acting by conjugation. The local orbital integral relevant to this relative trace formula is
\begin{equation}\label{orb int GL}
J^{\GL}_{x,\g}(f):=\int_{\GL_{n-1}(F_{x})}f(h^{-1}\g h)\eta_{x}(\det h)dh,\quad \g\in \Sg_{n}(F_{x}), f\in C_{c}^{\infty}(\Sg_{n}(F_{x})).
\end{equation}
Here $\eta_{x}$ is the character $F^{\times}_{x}\to \{\pm1\}$ attached to the quadratic extension $F'_{x}/F_{x}$. 

The Lie algebra analogue of $\GL_{n-1}\bs \Sg_{n}$ is $\GL_{n-1}\bs(\gl_{n}\otimes_{F}F'_{-})$ where $F'_{-}=(F')^{\s=-1}$, and $\GL_{n-1}$ acts by conjugation. Let $V_{n}$ be the standard representation of $\GL_{n}$ over $F$. It is more convenient to identify $\GL_{n-1}\bs (\gl_{n}\otimes_{F}F'_{-})$ with
\begin{equation*}
\GL_{n}\bs \left(\Hom_{F}(V_{n}, V_{n}\ot F'_{-})\times (V_{n} \times V^{*}_{n})^{1}\right),
\end{equation*}
where $(V_{n} \times V^{*}_{n})^{1}$ consists of $(e,e^{*})\in V_{n}\times V_{n}^{*}$ such that $e^{*}(e)=1$, and $\GL_{n}$ is acting diagonally on all factors (conjugation on the first factor).  The GIT quotient $\frC$ of $\Hom_{F}(V_{n}, V_{n}\ot F'_{-})\times (V_{n} \times V^{*}_{n})^{1}$ by $\GL_{n}$ is an affine space of dimension $2n-1$. For $(\ph, e, e^{*})\in \Hom_{F}(V_{n}, V_{n}\ot F'_{-})\times (V_{n} \times V^{*}_{n})^{1}$, we have invariants $a_{i}(\ph)\in (F'_{-})^{\ot i}$ that records the $i$-th coefficient of the characteristic polynomial of $\ph$ ($1\le i\le n$), and $b_{i}=e^{*}(\ph^{i}e)\in (F'_{-})^{\ot i}$ for $1\le i\le n-1$. The invariants $(a_{1},\cdots, a_{n},b_{1},\cdots, b_{n-1})$ give coordinates for $\frC$.

We introduce the following moduli stack $\cM$ which serves as a global avatar for the Lie algebra version of the orbital integrals appearing in this relative trace formula.  Fix line bundles $\cL$ and $\cL'$ on  $X$. Let $\cL_{-}=\cL\otimes_{\cO_{X}}\cO_{X'}^{\s=-1}$. The stack $\cM$ classifies tuples $(\cE, \ph, s,s^{*})$ where $\cE$ is a vector bundle of rank $n$ over $X$, $\ph: \cE\to \cE\otimes\cL_{-}$, $s: \cL'^{-1}\to \cE$ and $s^{*}: \cE\to \cL'$ are $\cO_{X}$-linear maps of coherent sheaves. The ``Hitchin base'' $\cB$ in this situation is the affine space $\prod_{i=1}^{n}\Gamma(X, \cL_{-}^{\ot i})\times\prod_{i=0}^{n-1}\Gamma(X, \cL'^{\ot 2}\ot \cL^{\ot i}_{-})$. The Hitchin map $f: \cM\to \cB$ sends $(\cE, \ph, s,s^{*})$ to the point of $\cB$ with coordinates $(a_{1}(\ph), \cdots, a_{n}(\ph), b_{0},\cdots, b_{n-1})$, where $a_{i}(\ph)$ are the coefficients of the characteristic polynomial of $\ph$, and $b_{i}=s^{*}\circ \ph^{i}\circ s: \cL'^{-1}\to \cL^{\ot i}_{-}\ot \cL'$.
\end{exam}

\begin{exam}\label{ex:JR Ug} Let $F'/F$ and $\nu: X'\to X$ be as in Example \ref{ex:JR GL}. Let $W_{n-1}$ be a Hermitian vector space of dimension $n-1$ over $F'$. Let $W_{n}=W_{n-1}\oplus F'e_{n}$ with the Hermitian form $(\cdot,\cdot)$ extending that on $W_{n-1}$ and such that  $W_{n-1}\bot e_{n}$, $(e_{n},e_{n})=1$.  Let $\Ug_{n}$ and $\Ug_{n-1}$ be the unitary groups over $F$ attached to $W_{n}$ and $W_{n-1}$. Consider $G'=\Ug_{n}\times \Ug_{n-1}$, and the subgroup $H'_{1}=H'_{2}=\Ug_{n-1}$  diagonally embedded into $G'$. 

The  double quotient $H'_{1}\bs G'/H'_{2}$ can be identified with the quotient  $\Ug_{n-1}\bs \Ug_{n}$ where $\Ug_{n-1}$ acts by conjugation. For $x\in |X|$, the local orbital integral relevant to this relative trace formula is
\begin{equation*}
J^{\Ug}_{x,\d}(f)=\int_{\Ug_{n-1}(F_{x})}f(h^{-1}\d h)dh, \quad \d\in \Ug_{n}(F_{x}), f\in C_{c}^{\infty}(\Ug_{n}(F_{x})).
\end{equation*}

The Lie algebra analogue of $\Ug_{n-1}\bs \Ug_{n}$ is $\Ug_{n-1}\bs \fru_{n}$, where $\fru_{n}$, the Lie algebra of $\Ug_{n}$, consists of skew-self-adjoint endomorphisms of $W_{n}$. As in the case of Example \ref{ex:JR GL}, we identify $\Ug_{n-1}\bs \fru_{n}$ with $\Ug_{n}\bs \left(\fru_{n}\times W^{1}_{n}\right)$ where $W_{n}^{1}$ is the set of vectors $e\in W_{n}$ such that $(e,e)=1$. The GIT quotient of $\fru_{n}\times W^{1}_{n}$ by $\Ug_{n}$ can be identified with the space $\frC$ introduced in Example \ref{ex:JR GL}. For $(\psi, e)\in \fru_{n}\times W^{1}_{n}$, its image in $\frC$ is $(a_{1}(\psi), \cdots, a_{n}(\psi), b_{1},\cdots, b_{n-1})$ where $a_{i}(\psi)\in (F'_{-})^{\ot i}$ are the coefficients of the characteristic polynomial of $\psi$ (as an endomorphism of $W_{n}$), and $b_{i}=(\psi^{i}e, e)\in (F'_{-})^{\ot i}$, since $\s(\psi^{i}e, e)=(e,\psi^{i}e)=(-1)^{i}(\psi^{i}e,e)$.

We introduce a moduli stack $\cN$ which serves as a global avatar for the Lie algebra version of the orbital integrals appearing in this relative trace formula. Fix line bundles $\cL$ and $\cL'$ on  $X$. The stack $\cN$ classifies tuples $(\cF, h, \psi, t)$ where $\cF$ is a vector bundle of rank $n$ on $X'$, $h: \cF\isom \s^{*}\cF^{\vee}$ is a Hermitian form on $\cF$, $\psi: \cF\to \cF\otimes \nu^{*}\cL$  is skew-self-adjoint with respect to $h$ and $t: \nu^{*}\cL'^{-1}\to \cF$ is an $\cO_{X'}$-linear map. When $\nu$ is unramified, the base $\cB$ introduced in Example \ref{ex:JR GL} still serves as the Hitchin base for $\cN$. The Hitchin map $g: \cN\to \cB$ sends $(\cF, h, \psi, t)$ to $(a_{1}(\psi), \cdots, a_{n}(\psi), b_{0},\cdots, b_{n-1})$, where $a_{i}(\psi)$ are the coefficients of the characteristic polynomial of $\psi$, and $b_{i}= t^{\vee}\circ h\circ\psi^{i}\circ t$ descends to an $\cO_{X}$-linear map  $\cL'^{-1}\to \cL^{\ot i}_{-}\ot \cL'$.
\end{exam}

\begin{theorem}[\cite{YJR}]\label{th:JR} Let $x$ be a place of $F$ such that $F'/F$ is unramified over $x$ and the Hermitian space $W_{n,x}$ has a self-dual lattices $\L_{n,x}$. Then for strongly regular semisimple elements $\g\in \Sg_{n}(F_{x})$ and $\d\in \Ug_{n}(F_{x})$ with the same invariants in $\frC(F_{x})$, we have
\begin{equation*}
J^{\GL}_{x,\g}(\one_{\Sg_{n}(\cO_{x})})=\pm J^{\Ug}_{x,\d}(\one_{\Ug(\L_{n,x})})
\end{equation*}
for some sign depending on the invariants of $\g$.
\end{theorem}

The main geometric observation in \cite{YJR} is that both $f:\cM\to \cB$ and $g:\cN\to \cB$ are small maps when restricted to a certain open subset of $\cB$. This enables us to prove an isomorphism between the direct images complexes of $f$ and $g$ by checking over the generic point of $\cB$. Such an isomorphism of sheaves, after passing to the Frobenius traces on stalks,  implies the identity above, which was the fundamental lemma conjectured by Jacquet and Rallis.

\section{Hitchin moduli stack and Shtukas}\label{s:Sht}

In this section we consider automorphic objects that arise as cohomology classes of moduli stacks of Shtukas, which are the function-field counterpart of Shimura varieties. These cohomology classes generalize the notion of automorphic forms. The periods and relative traces have their natural analogues in this more general setting. Hitchin-type moduli stacks continue to play an important role in the study of such relative trace formulae. We give a survey of our recent work \cite{YZ1},\cite{YZ2} on higher Waldspurger-Gross-Zagier formulae and \cite{YAFL} on the arithmetic fundamental lemma, which fit into the framework to be discussed in this section.

\subsection{Moduli of Shtukas}\label{ss:Sht}
In his seminal paper \cite{DrEll}, Drinfeld introduced the moduli of elliptic modules as a function field analogue of modular curves. Later in \cite{DrSht}, Drinfeld defined more general geometric object called Shtukas, and used them to prove the Langlands conjecture for $\GL_{2}$ over function fields. Since then it became clear that the moduli stack of Shtukas should play the role of Shimura varieties for function fields, and its cohomology should realize the Langlands  correspondence for global function fields. This idea was realized for $\GL_{n}$ by L.Lafforgue \cite{LLaff} who proved the full Langlands conjecture in this case. For an arbitrary reductive group $G$, V.Lafforgue \cite{VLaff} proved the automorphic to Galois direction of the Langlands conjecture using moduli stacks of Shtukas.

\sss{The moduli of Shtukas} The general definition of $G$-Shtukas was given by Varshavsky \cite{Var}. For simplicity of presentation we assume $G$ is split. Again we fix an open subgroup $K\subset K_{0}$, and let $N\subset |X|$ be the finite set of places where $K_{x}\ne G(\cO_{x})$. Choosing a maximal split torus $T$ and a Borel subgroup $B$ containing $T$, we may therefore talk about dominant coweights of $T$ with respect to $B$.  Let $r\ge0$ be an integer. Let $\mu=(\mu_{1},\cdots, \mu_{r})$ be a sequence of dominant coweights of $T$. Recall dominant coweights of $T$ are in bijection with relative positions of two $G$-bundles over the formal disk with the same generic fiber. 

Let $\Hk^{\mu}_{G,K}$ be the Hecke stack classifying points $x_{1},\cdots, x_{r}\in X-N$ together with a diagram of the form
\begin{equation*}
\xymatrix{\cE_{0}\ar@{-->}[r]^{f_{1}} & \cE_{1}\ar@{-->}[r]^{f_{2}} & \cdots \ar@{-->}[r]^{f_{r}} & \cE_{r}}
\end{equation*}
where $\cE_{i}$ are $G$-bundles over $X$ with $K$-level structures, and $f_{i}: \cE_{i-1}|_{X-x_{i}}\isom \cE_{i}|_{X-x_{i}}$ is an isomorphism compatible with the level structures whose relative position at $x_{i}$ is in the closure of that given by $\mu_{i}$.

A {\em $G$-Shtuka of type $\mu$ with level $K$} is the same data as those classified by $\Hk^{\mu}_{G, K}$, together with an isomorphism of $G$-bundles compatible with $K$-level structures
\begin{equation}\label{ErE0}
\io: \cE_{r}\isom {}^{\tau}\cE_{0}.
\end{equation}
Here, ${}^{\tau}\cE_{0}$ is the image of $\cE_{0}$ under the Frobenius morphism $\Fr: \Bun_{G,K}\to \Bun_{G,K}$. If we are talking about an $S$-family of $G$-Shtukas for some $k$-scheme $S$, $\cE_{0}$ is a $G$-torsor over $X\times S$, then ${}^{\tau}\cE_{0}:=(\id_{X}\times \Fr_{S})^{*}\cE_{0}$. There is a moduli stack $\Sht^{\mu}_{G,K}$ of $G$-Shtuka of type $\mu$, which fits into a Cartesian diagram
\begin{equation}\label{defn ShtG}
\xymatrix{\Sht^{\mu}_{G,K}\ar[r]\ar[d] & \Hk^{\mu}_{G,K}\ar[d]^{(p_{0}, p_{r})}\\
\Bun_{G,K} \ar[r]^-{(\id, \Fr)}& \Bun_{G,K}\times \Bun_{G,K}
}
\end{equation}
Let us observe the similarity with the definition of the (group version of) Hitchin stack in diagram \eqref{MHk geom}: the main difference is that we are replacing the diagonal map of $\Bun_{G,K}$ by the graph of the Frobenius.

Recording only the points $x_{1},\cdots, x_{r}$ gives a morphism
\begin{equation*}
\pi^{\mu}_{G,K}: \Sht^{\mu}_{G,K}\to (X-N)^{r}
\end{equation*}

The datum $\mu$ is called {\em admissible} if $\sum_{i}\mu_{i}$ lies in the coroot lattice. The existence of an isomorphism \eqref{ErE0} forces $\mu$ to be admissible. Therefore $\Sht^{\mu}_{G,K}$ is nonempty only when $\mu$ is admissible.  When $r=0$, $\Sht^{\mu}_{G,K}$ is the discrete stack given by the double coset $\Bun_{G,K}(k)=G(F)\bs G(\AA)/K$. For $\mu$ admissible, we have
\begin{equation*}
d_{G}(\mu):=\dim \Sht^{\mu}_{G,K}=\sum_{i}(\jiao{2\r_{G}, \mu_{i}}+1).
\end{equation*}

\sss{Hecke symmetry}\label{sss:Sht Hk}
Let $g\in G(\AA)$ and let $S$ be the finite set of $x\in |X|-N$ such that $g_{x}\notin G(\cO_{x})$. There is a self-correspondence $\Sht^{\mu}_{G,KgK}$ (the dependence on $g$ is only through the double coset $KgK$) of $\Sht^{\mu}_{G,K}|_{(X-N-S)^{r}}$ such that both maps to $\Sht^{\mu}_{G,K}|_{(X-N-S)^{r}}$ are finite \'etale. It then induces an endomorphism of the direct image complex $\bR\pi^{\mu}_{G,K,!}\IC(\Sht^{\mu}_{G,K})|_{(X-N-S)^{r}}$.  V.Lafforgue \cite{VLaff} used his construction of excursion operators to extend this endomorphism to the whole complex $\bR\pi^{\mu}_{G,K,!}\IC(\Sht^{\mu}_{G,K})$ over $(X-N)^{r}$. If we assign this endomorphism to the function $\one_{KgK}$, it extends by linearity to an action of the Hecke algebra $C_{c}(K\bs G(\AA)/K)$ on the complex $\bR\pi^{\mu}_{G,K,!}\IC(\Sht^{\mu}_{G,K})$, and hence on its geometric stalks and on its cohomology groups.

\sss{Intersection cohomology of $\Sht^{\mu}_{G,K}$} The singularities of the map $\pi^{\mu}_{G,K}$ are exactly the same as the product of the Schubert varieties $\Gr_{G,\le\mu_{i}}$ in the affine Grassmannian $\Gr_{G}$. It is expected that the complex $\bR\pi^{\mu}_{G,K,!}\IC(\Sht^{\mu}_{G,K})$ should realize the global Langlands correspondence for $G$ in a way similar to the Eichler-Shimura correspondence for modular curves. The phenomenon of endoscopy makes stating a precise conjecture quite subtle, but a rough form of the expectation is a $C_{c}(K\bs G(\AA)/K)$-equivariant decomposition over $(X-N)^{r}$
\begin{equation}\label{coho Sht}
\bR\pi^{\mu}_{G,K,!}\IC(\Sht^{\mu}_{G, K})''=''\left(\bigoplus_{\pi \textup{ cuspidal}} \pi^{K}\otimes (\boxtimes_{i=1}^{r}\rho^{\mu_{i}}_{\pi}[1])\right)\bigoplus (\textup{Eisenstein part}).
\end{equation}
Here $\pi$ runs over cuspidal automorphic representations of $G(\AA)$ such that $\pi^{K}\ne0$, $\rho_{\pi}$ is the $\dG$-local system on $X-N$ attached to $\pi$ by the Langlands correspondence, and $\rho^{\mu_{i}}_{\pi}$ is the local system obtained by the composition
\begin{equation*}
\pi_{1}(X-N, *)\xr{\rho_{\pi}}\dG(\Qlbar)\to \GL(V(\mu_{i}))
\end{equation*}
where $V(\mu_{i})$ is the irreducible representation of the dual group $\dG$ with highest weight $\mu_{i}$. 

One approach to prove \eqref{coho Sht} is to use trace formulae. One the one hand, consider the action of a Hecke operator composed with a power of Frobenius at some $x\in |X|-N$ acting on the geometric stalk at $x$ of the left side of \eqref{coho Sht}, which is $\IHc{*}{\Sht^{\mu}_{G,K,\ov x}}$. The trace of this action can be calculated by the Lefschetz trace formula, and can be expressed as a sum of twisted orbital integrals. On the other hand, the trace of the same operator on the right side of \eqref{coho Sht} can be calculated by the Arthur-Selberg trace formula, and be expressed using orbital integrals. The identity \eqref{coho Sht} would then follow from an identity between the twisted orbital integrals and the usual orbital integrals that appear in both trace formulae, known as the base-change fundamental lemma.

The difficulty in implementing this strategy is that $\Sht^{\mu}_{G,K}$ is not of finite type, and both the Lefschetz trace and the Arthur-Selberg trace would be divergent. L.Lafforgue \cite{LLaff} treated the case $G=\GL_{n}$ and $\mu=((1,0,\cdots,0), (0,\cdots,0,-1))$ by difficult analysis of the compactification of truncations of $\Sht^{\mu}_{G,K}$, generalizing the work of Drinfeld on $\GL_{2}$.

\sss{Cohomological spectral decomposition} We discuss a weaker version of the spectral decomposition \eqref{coho Sht}.  As mentioned above, $\Sht^{\mu}_{G,K}$ is not of finite type, so its intersection cohomology is not necessarily finite-dimensional.  One can present $\Sht^{\mu}_{G,K}$ as an increasing union of finite-type open substacks, but these substacks are not preserved by the Hecke correspondences. Despite all that, we expect nice finiteness properties of $\IHc{*}{\Sht^{\mu}_{G,K}\ot\kbar}$ as a Hecke module. More precisely,   the spherical Hecke algebra $C_{c}(K^{N}_{0}\bs G(\AA^{N})/K_{0}^{N},\Qlbar)$ (superscript $N$ means removing places in $N$) should act through a quotient $\ov\sH^{N}$ (possibly depending on $\mu$) which is a finitely generated algebra over $\Qlbar$, and that $\IHc{*}{\Sht^{\mu}_{G,K}\ot\kbar}$ should be a finitely generated module over $\ov\sH^{N}$. Now viewing  $\IHc{*}{\Sht^{\mu}_{G,K}\ot\kbar}$ as a coherent sheaf on $\ov\sH^{N}$, we should get a canonical decomposition of it in terms of connected components of $\ov\sH^{N}$. A coarser decomposition should take the following form
\begin{equation}\label{coho spec}
\IHc{*}{\Sht^{\mu}_{G,K}\ot\kbar}=\bigoplus_{[P]}\IHc{*}{\Sht^{\mu}_{G,K}\ot\kbar}_{[P]}
\end{equation}
where $[P]$ runs over associated classes of parabolic subgroups of $G$. The support of $\IHc{*}{\Sht^{\mu}_{G,K}\ot\kbar}_{[P]}$ should be described using the analogous quotient of the Hecke algebra $\ov\sH^{N}_{L}$ for the Levi factor $L$ of $P$, via the Satake transform from the spherical Hecke algebra for $G$ to the one for  $L$.  If $G$ is semisimple, the part $\IHc{*}{\Sht^{\mu}_{G,K}\ot\kbar, \Qlbar}_{[G]}$ should be finite-dimensional over $\Qlbar$.  

In the simplest nontrivial case $G=\PGL_{2}$ we have proved the coarse decomposition.

\begin{theorem}[\cite{YZ1},\cite{YZ2}]\label{th:coho spec} For $G=\PGL_{2}$, consider the moduli of Shtukas $\Sht^{r}_{G}$ without level structures of type $\mu=(\mu_{1},\cdots, \mu_{r})$ where each $\mu_{i}$ is the minuscule coweight. Then there is a decomposition of Hecke modules
\begin{equation*}
\cohoc{2r}{\Sht^{r}_{G}\ot\kbar}=(\oplus_{\chi} \cohoc{2r}{\Sht^{r}_{G}\ot\kbar}[\chi])\oplus \cohoc{2r}{\Sht^{r}_{G}\ot\kbar}_{\Eis}
\end{equation*}
where $\chi$ runs over a finite set of characters of the Hecke algebra $C_{c}(K_{0}\bs G(\AA)/K_{0})$, and the support of $\cohoc{2r}{\Sht^{r}_{G}\ot\kbar}_{\Eis}$ is defined by the Eisenstein ideal. 

For $i\ne 2r$, $\cohoc{i}{\Sht^{r}_{G}\ot\kbar}$ is finite-dimensional.

Similar result holds for a version of $\Sht^{r}_{G}$ with Iwahori level structures.
\end{theorem}
We expect the similar techniques to work for general split $G$ and general type $\mu$. 

Assume we have an analogue of the above theorem for $G$. Let $\pi$ be a cuspidal automorphic representation of $G(\AA)$ such that $\pi^{K}\ne0$, then  $C_{c}(K^{N}_{0}\bs G(\AA^{N})/K_{0}^{N},\Qlbar)$ acts on $\pi^{K}$ by a character $\chi_{\pi}$ up to semisimplification. Suppose $\chi_{\pi}$ does not appear in the support of $\IHc{*}{\Sht^{\mu}_{G}\ot\kbar}_{[P]}$ for any proper parabolic $P$ (in which case we say $\chi_{\pi}$ is non-Eisenstein), then the generalized eigenspace  $\IHc{*}{\Sht^{\mu}_{G,K}\ot\kbar}[\chi_{\pi}]$ is a finite-dimensional direct summand of $\IHc{*}{\Sht^{\mu}_{G,K}\ot\kbar}$ containing the contribution of $\pi$ but possibly also companions of $\pi$ with the same Hecke character $\chi_{\pi}$ away from $N$.

\subsection{Heegner-Drinfeld cycles and periods}\label{ss:HD}
When $r=0$, the left side of \eqref{coho Sht} is simply the function space $C_{c}(G(F)\bs G(\AA)/K)$ where cuspidal automorphic forms live.   In general, we should think of  cohomology classes in $\IHc{*}{\Sht^{\mu}_{G,K,\ov x}}$ or $\IHc{*}{\Sht^{\mu}_{G,K}\ot\kbar}$ as generalizations of automorphic forms.  We shall use this viewpoint to generalize some constructions in Section \ref{s:TF} from classical automorphic forms to cohomology classes of $\Sht^{\mu}_{G,K}$.

\sss{Heegner-Drinfeld cycles}  Let $H\subset G$ be a subgroup defined over $k$ with level group $K_{H}=K\cap H(\AA)$. It induces a map $\th_{\Bun}:\Bun_{H,K_{H}}\to \Bun_{G,K}$. 

Fix an integer $r\ge0$. Let $\l=(\l_{1},\cdots, \l_{r})$ be an admissible sequence of dominant coweights of $H$; let $\mu=(\mu_{1},\cdots,\mu_{r})$ be an admissible sequence of dominant coweights of $G$. To relate $\Sht^{\l}_{H}$ to $\Sht^{\mu}_{G}$, we need to impose more restrictions on $\l$ and $\mu$. For two coweights $\mu,\mu'$ of $G$ we write $\mu\le_{G}\mu'$ if for some (equivalently all) choices of a Borel $B'\subset G_{\kbar}$ and a maximal torus $T'\subset B'$, $\mu'_{B',T'}-\mu_{B',T'}$ is a sum of positive roots. Here $\mu_{B',T'}$ (resp. $\mu'_{B',T'}$) is the unique dominant coweight of $T'$ conjugate to $\mu$ (resp. $\mu'$). 

We assume that $\l_{i}\le_{G} \mu_{i}$ for  $0\le i \le r$. In this case there is a natural morphism of Hecke stacks
\begin{equation*}
\th_{\Hk}: \Hk^{\l}_{H,K_{H}}\to \Hk^{\mu}_{G,K}
\end{equation*}
compatible with $\th_{\Bun}$. The Cartesian diagram \eqref{defn ShtG} and its counterpart for $\Sht^{\l}_{H}$ then induce a map over $(X-N)^{r}$
\begin{equation*}
\th: \Sht^{\l}_{H,K_{H}}\to \Sht^{\mu}_{G,K}.
\end{equation*}
If $\th$ is proper, the image of the fundamental class of $\Sht^{\l}_{H,K_{H}}$ defines an algebraic cycle $\Sht^{\mu}_{G,K}$ which we call a {\em Heegner-Drinfeld cycle}.

\begin{exam}\label{ex:T} Consider the case $G=\PGL_{2}$, and $H=T$ is a non-split torus of the form $T=(\Res_{F'/F}\Gm)/\Gm$ for some quadratic extension $F'/F$. Since $T$ is not a constant group scheme over $X$, our previous discussion does not directly apply, but we can easily define what a $T$-Shtuka is. The quadratic extension $F'$ is the function field of a smooth projective curve $X'$ with a degree two map $\nu: X'\to X$. Let $\l=(\l_{1},\cdots, \l_{r})\in \ZZ^{r}$ with $\sum_{i}\l_{i}=0$, we may consider the moduli of rank one Shtukas $\Sht^{\l}_{\GL_{1}, X'}$ over $X'$ of type $\l$. We define $\Sht^{\l}_{T}$ to be the quotient $\Sht^{\l}_{\GL_{1}, X'}/\Pic_{X}(k)$, where the discrete groupoid $\Pic_{X}(k)$ is acting by pulling back to $X'$ and tensoring with rank one Shtukas. It can be shown that the projection $\Sht^{\l}_{T}\to X'^{r}$ is a finite \'etale Galois cover with Galois group $\Pic_{X'}(k)/\Pic_{X}(k)$. In particular, $\Sht^{\l}_{T}$ is a smooth and proper DM stack over $k$ of dimension $r$.  

Now let $\mu=(\mu_{1},\cdots, \mu_{r})$ be a sequence of dominant coweights of $G$. Then each $\mu_{i}$ can be identified with an element in $\ZZ_{\ge0}$, with the positive coroot corresponding to $1$.   Admissibility of $\mu$ means that $\sum_{i}\mu_{i}$ is even. The condition $\l_{i}\le_{G} \mu_{i}$  is saying that $|\l_{i}|\le \mu_{i}$ and that $\mu_{i}-\l_{i}$ is even. When $\l_{i}\le_{G}\mu_{i}$ for all $i$, the map $\th: \Sht^{\l}_{T}\to \Sht^{\mu}_{G}$ simply takes a rank one Shtuka $(\{\cE_{i}\};\{x'_{i}\})$ on $X'$ and sends it to the direct image $(\{\nu_{*}\cE_{i}\};\{\nu(x'_{i})\})$, which is a rank two Shtuka on $X$.
\end{exam}

\sss{Periods}\label{sss:Sht per}
Fix a Haar measure $\mu_{H}$ on $H(\AA)$. Under a purely root-theoretic condition on $\l$ and $\mu$,  $\th^{*}$ induces a map
\begin{equation*}
\th^{*}: \IHc{2d_{H}(\l)}{\Sht^{\mu}_{G,K}\ot\kbar}\to\cohoc{2d_{H}(\l)}{\Sht^{\l}_{H,K_{H}}\ot\kbar},
\end{equation*}
and therefore defines a period map
\begin{equation*}
\cP^{G,\mu}_{H,\l}: \IHc{2d_{H}(\l)}{\Sht^{\mu}_{G,K}\ot\kbar}\xr{\th^{*}}\cohoc{2d_{H}(\l)}{\Sht^{\l}_{H,K_{H}}\ot\kbar}\xr{\cap[\Sht^{\l}_{H,K_{H}}]\cdot\vol(K_{H},\mu_{H})}\Qlbar.
\end{equation*}
The last map above is the cap product with the fundamental class of $\Sht^{\l}_{H}$ followed by multiplication by $\vol(K_{H},\mu_{H})$.

Now assume $\Sht^{\l}_{H,K_{H}}$ has half the dimension of $\Sht^{\mu}_{G,K}$, 
\begin{equation}
d_{G}(\mu)=2d_{H}(\l).
\end{equation}
Let $\pi$ be a cuspidal automorphic representation of $G(\AA)$. To make sense of periods on $\pi$, we assume for the moment that the contribution of $\pi$ to the intersection cohomology of $\Sht^{\mu}_{G,K}$ is as predicted in \eqref{coho Sht}. Restricting $\cP^{G,\mu}_{H,\l}$ to the $\pi$-part we get
\begin{equation*}
\pi^{K}\ot(\otimes_{i=1}^{r}\cohoc{1}{(X-N)\ot\kbar, \rho^{\mu_{i}}_{\pi}})\xr{\th^{*}}\cohoc{2d_{H}(\l)}{\Sht^{\l}_{H,K_{H}}\ot\kbar}\to\Qlbar.
\end{equation*} 
As we expect $\rho^{\mu_{i}}_{\pi}$ to be pure, the above map should factor through the pure quotient of $\cohoc{1}{(X-N)\ot\kbar, \rho^{\mu_{i}}_{\pi}}$, which is $\cohoc{1}{X\ot{\kbar}, j_{!*}\rho^{\mu_{i}}_{\pi}}$ (the cohomology of the middle extension of $\rho^{\mu_{i}}_{\pi}$), and which does not change after enlarging $N$. Now shrinking $K$ and passing to the direct limit, we get 
\begin{equation*}
\cP^{G,\mu}_{H,\l, \pi}: \pi\ot(\otimes_{i=1}^{r}\cohog{1}{X\ot{\kbar},  j_{!*}\rho^{\mu_{i}}_{\pi}})\to \Qlbar
\end{equation*}
which is the analogue of the classical period \eqref{per}. Again $\cP^{G,\mu}_{H,\l, \pi}$ should factor through the coinvariants $\pi_{H(\AA)}\otimes(\cdots)$.


\subsection{Shtuka version of relative trace formula}\label{ss:Sht RTF} 
\sss{The setup} Let $H_{1}$ and $H_{2}$ be reductive subgroups of $G$ over $k$. Fix an integer $r\ge0$. Let $\l=(\l_{1},\cdots, \l_{r})$ (resp. $\k$ and $\mu$) be an admissible sequence of dominant coweights of $H_{1}$ (resp. $H_{2}$ and $G$). Assume that $\l_{i}\le_{G} \mu_{i}$ and $\k_{i}\le_{G}\mu_{i}$. In this case there are natural morphisms 
\begin{equation*}
\xymatrix{\Sht^{\l}_{H_{1}}\ar[r]^{\th_{1}} & \Sht^{\mu}_{G} & \Sht^{\k}_{H_{2}}\ar[l]_{\th_{2}}}
\end{equation*}
Suppose
\begin{equation*}
\dim \Sht^{\l}_{H_{1}}=\dim \Sht^{\k}_{H_{2}}=\frac{1}{2}\dim \Sht^{\mu}_{G}.
\end{equation*}
With the same extra assumptions as in Section \ref{sss:Sht per}, we may define the periods $\cP^{G,\mu}_{H_{1},\l, \pi}$ and $\cP^{G,\mu}_{H_{2},\k, \wt\pi}$ (where $\wt\pi$ is the contragradient of $\pi$). We expect the tensor product 
\begin{equation*}
\cP^{G,\mu}_{H_{1},\l, \pi}\ot\cP^{G,\mu}_{H_{2},\k, \wt\pi}: \pi\ot\wt\pi\ot\left(\ot_{i=1}^{r}\left(\cohog{1}{X\ot{\kbar},  j_{!*}\rho^{\mu_{i}}_{\pi}}\ot\cohog{1}{X\ot{\kbar},  j_{!*}\rho^{\mu_{i}}_{\wt\pi}}\right)\right)\to \Qlbar
\end{equation*}
to factor through the pairing between $\cohog{1}{X\ot{\kbar},  j_{!*}\rho^{\mu_{i}}_{\pi}}$ and $\cohog{1}{X\ot{\kbar},  j_{!*}\rho^{\mu_{i}}_{\wt\pi}}$ given by the cup product (the local systems $\rho^{\mu_{i}}_{\wt\pi}$ and $\rho^{\mu_{i}}_{\pi}$ are dual to each other up to a Tate twist).  Assuming this, we get a pairing 
\begin{equation*}
\cP^{G,\mu}_{H_{1},\l, \pi}\ot\cP^{G,\mu}_{H_{2},\k, \wt\pi}: \pi_{H_{1}(\AA)}\ot\wt\pi_{H_{2}(\AA)}\to \Qlbar.
\end{equation*}
When $\dim\pi_{H_{1}(\AA)}=\dim\wt\pi_{H_{2}(\AA)}=1$, we expect the ratio between the above pairing and the Petersson inner product to be related to {\em derivatives} of $L$-functions of $\pi$, though I do not know how to formulate a precise conjecture in general.

\sss{The relative trace} One way to access  $\cP^{G,\mu}_{H_{1},\l, \pi}\ot\cP^{G,\mu}_{H_{2},\k, \wt\pi}$ is to develop a relative trace formula whose spectral expansion gives these periods. Fix a compact open subgroup $K\subset G(\AA)$ and let $K_{i}=H_{i}(\AA)\cap K$. Assume $\th_{i}$ are proper, we have Heegner-Drinfeld cycles
\begin{equation*}
\cZ^{\l}_{H_{1}}=\th_{1*}[\Sht^{\l}_{H_{1},K_{1}}], \quad \cZ^{\k}_{H_{2}}=\th_{2*}[\Sht^{\k}_{H_{2},K_{2}}]
\end{equation*}
both of half dimension in $\Sht^{\mu}_{G,K}$. Intuitively we would like to form the intersection number
\begin{equation}\label{defn I}
I^{(G,\mu)}_{(H_1,\l),(H_2,\k)}(f)=\jiao{\cZ^{\l}_{H_1}, f*\cZ^{\k}_{H_2}}_{\Sht^{\mu}_{G,K}}, \quad f\in C_{c}(K\bs G(\AA)/K)
\end{equation}
as the ``relative trace'' of $f$ in this context. Here $f*(-)$ denotes the action of the Hecke algebra on the Chow group of $\Sht^{\mu}_{G,K}$, defined similarly as in Section \ref{sss:Sht Hk}. However, there are several technical issues before we can make sense of this intersection number.
\begin{enumerate}
\item $\Sht^{\mu}_{G,K}$ may not be smooth so the intersection product of cycles may not be defined. 
\item Suppose the intersection of $\cZ^{\l}_{H_1}$ and $\cZ^{\k}_{H_2}$ is defined as a $0$-cycle on $\Sht^{\mu}_{G,K}$, if we want to get a number out of this $0$-cycle, we need it to be a proper cycle, i.e., it should lie in the Chow group of cycles with proper support (over $k$).
\end{enumerate}
The first issue goes away if we assume each $\mu_{i}$ to be a minuscule coweight of $G$, which guarantees that $\Sht^{\mu}_{G,K}$ is smooth over $(X-N)^{r}$. The second issue is more serious and is analogous to the divergence issue for the usual relative trace. In results that we will present later, it won't be an issue because there $\cZ^{\l}_{H_1,K_{1}}$ is itself a proper cycle. In the sequel we will proceed with heuristic arguments as we did in Section \ref{ss:RTF}, and ignore these issues.

%

When $\l,\k,\mu$ are all zero, the linear functional $I^{(G,0)}_{(H_{1},0),(H_{2},0)}$ becomes the relative trace $\RTr^{G}_{H_{1}, H_{2}}$ defined in \eqref{rel trace}. Therefore the functional $I^{(G,\mu)}_{(H_1,\l),(H_2,\k)}$ is a generalization of the relative trace.

\sss{Intersection number in terms of Hitchin-like moduli stacks} In the case of the usual relative trace for $f=\one_{KgK}$, we introduced a Hitchin-like moduli stack $\cM^{G}_{H_{1}, H_{2},KgK}$ whose point-counting is essentially the relative trace of $f$. We now try to do the same for $I^{(G,\mu)}_{(H_1,\l),(H_2,\k)}$. To simplify notations we assume $K=K_{0}=\prod_{x} G(\cO_{x})$, hence $K_{i}=\prod H_{i}(\cO_{x})$, and suppress them from the notation for Shtukas.

To calculate the intersection number \eqref{defn I}, a natural starting point is to form the stack-theoretic intersection of the cycles $\cZ^{\l}_{H_1}$ and $f*\cZ^{\k}_{H_2}$, i.e., consider the Cartesian diagram
\begin{equation}\label{defn ShtM}
\xymatrix{\Sht^{(G,\mu)}_{(H_1,\l),(H_2,\k),K_{0}gK_{0}} \ar[rr]\ar[d] && \Sht^{\mu}_{G,K_{0}gK_{0}}\ar[d]\\
\Sht^{\l}_{H_1}\times\Sht^{\k}_{H_2}\ar[rr]^-{\th_{1}\times\th_{2}} && \Sht^{\mu}_{G}\times\Sht^{\mu}_{G}}
\end{equation}
The expected dimension of $\Sht^{(G,\mu)}_{(H_1,\l),(H_2,\k),K_{0}gK_{0}}$ is zero. If $\Sht^{(G,\mu)}_{(H_1,\l),(H_2,\k),K_{0}gK_{0}}$ indeed was zero-dimensional and moreover was proper over $k$, then $I^{(G,\mu)}_{(H_1,\l),(H_2,\k)}(\one_{K_{0}gK_{0}})$ would be equal to the length of $\Sht^{(G,\mu)}_{(H_1,\l),(H_2,\k),K_{0}gK_{0}}$. However, neither the zero-dimensionality nor the properness is true in general. Putting these issues aside, we proceed to rewrite $\Sht^{(G,\mu)}_{(H_1,\l),(H_2,\k),K_{0}gK_{0}}$ in Hitchin-like terms.

Recall we have Hecke correspondences $\Hk^{\mu}_{G}, \Hk^{\l}_{H_1}$ and $\Hk^{\k}_{H_2}$ for $\Bun_{G}, \Bun_{H_{1}}$ and $\Bun_{H_{2}}$ related by the maps $\th_{1,\Hk}$ and $\th_{2,\Hk}$. We can also define a Hecke correspondence for $\Hk_{G,K_{0}gK_{0}}$ as the moduli stack classifying $x_{1},\cdots, x_{r}\in X$ and a commutative diagram of rational isomorphisms of $G$-bundles over $X$ 
\begin{equation}\label{ladder}
\xymatrix{ \cE_{0}\ar@{-->}[r]\ar@{-->}[d]^{\ph_{0}} & \cE_{1}\ar@{-->}[r]\ar@{-->}[d]^{\ph_{1}} & \cdots\ar@{-->}[r] & \cE_{r}\ar@{-->}[d]^{\ph_{r}}\\
\cE'_{0}\ar@{-->}[r] & \cE'_{1}\ar@{-->}[r] & \cdots\ar@{-->}[r] & \cE'_{r}
}
\end{equation}
such that 
\begin{enumerate}
\item The top and bottom rows of the diagram give objects in $\Hk^{\mu}_{G}$ over $(x_{1},\cdots, x_{r})\in X^{r}$;
\item Each column of the diagram gives an object in $\Hk_{G, K_{0}gK_{0}}$, i.e., the relative position of $\ph_{i}$ is given by $K_{0}gK_{0}$ for $0\le i\le r$.
\end{enumerate}
We denote the resulting moduli stack by $\Hk^{\mu}_{\Hk,K_{0}gK_{0}}$. We have maps
$\frp, \frq : \Hk^{\mu}_{\Hk,K_{0}gK_{0}}\to \Hk^{\mu}_{G}$ by taking the top and the bottom rows; we also have maps $\frp_{i}: \Hk^{\mu}_{\Hk,K_{0}gK_{0}}\to \Hk_{G,K_{0}gK_{0}}$ by taking the $i$th column. The Hecke correspondence $\Sht^{\mu}_{G,K_{0}gK_{0}}$ in Section \ref{sss:Sht Hk} is defined as the pullback of $(\frp_{0},\frp_{r}):\Hk^{\mu}_{\Hk,K_{0}gK_{0}}\to \Hk_{G,K_{0}gK_{0}}\times\Hk_{G,K_{0}gK_{0}}$ along the graph of the Frobenius morphism for $\Hk_{G,K_{0}gK_{0}}$.

We then define $\Hk^{\l,\k}_{\cM, K_{0}gK_{0}}$ using the Cartesian diagram
\begin{equation*}
\xymatrix{\Hk^{\l,\k}_{\cM, K_{0}gK_{0}}\ar[rr]\ar[d] && \Hk^{\mu}_{\Hk,K_{0}gK_{0}}\ar[d]^{(\frp,\frq)}\\
\Hk^{\l}_{H_1}\times \Hk^{\k}_{H_2}\ar[rr]^-{\th_{1,\Hk}\times\th_{2,Hk}} && \Hk^{\mu}_{G}\times\Hk^{\mu}_{G}}
\end{equation*}
Now $\Hk^{\l,\k}_{\cM, K_{0}gK_{0}}$ can be viewed as an $r$-step Hecke correspondence for $\cM^{G}_{H_{1}, H_{2},K_{0}gK_{0}}$. Indeed, $\Hk^{\l,\k}_{\cM, K_{0}gK_{0}}$ classifies a diagram similar to \eqref{ladder}, except that $\cE_{i}$ (resp. $\cE_{i}'$) are now induced from $H_{1}$-bundles $\cF_{i}$ (resp. $H_{2}$-bundles $\cF'_{i}$), and the top row (resp. bottom row) are induced from an object in $\Hk^{\l}_{H_1}$ (resp. $\Hk^{\k}_{H_2}$). Now each column in such a diagram gives an object in $\cM^{G}_{H_{1}, H_{2},K_{0}gK_{0}}$. Recording the $i$-th column gives a map $m_{i}: \Hk^{\l,\k}_{\cM,K_{0}gK_{0}}\to \cM^{G}_{H_{1}, H_{2},K_{0}gK_{0}}$. We claim that there is a Cartesian diagram expressing $\Sht^{(G,\mu)}_{(H_1,\l),(H_2,\k),K_{0}gK_{0}}$ as the ``moduli of Shtukas for $\cM^{G}_{H_{1}, H_{2},K_{0}gK_{0}}$ with modification type $(\l,\k)$'':
\begin{equation}\label{alt ShtM}
\xymatrix{\Sht^{(G,\mu)}_{(H_1,\l),(H_2,\k),K_{0}gK_{0}} \ar[rr]\ar[d] && \Hk^{\l,\k}_{\cM,K_{0}gK_{0}}\ar[d]^{(m_{0},m_{r})}\\
\cM^{G}_{H_{1}, H_{2},K_{0}gK_{0}}\ar[rr]^-{(\id,\Fr)} && \cM^{G}_{H_{1}, H_{2},K_{0}gK_{0}}\times \cM^{G}_{H_{1}, H_{2},K_{0}gK_{0}}
}
\end{equation}
Indeed, this diagram is obtained by unfolding each corner of \eqref{defn ShtM} as a fiber product of an $r$-step Hecke correspondence with the graph of Frobenius, and re-arranging the order of taking fiber products.

Continuing with the heuristics, the intersection number $I^{(G,\mu)}_{(H_1,\l),(H_2,\k)}(\one_{K_{0}gK_{0}})$, which ``is'' the length of $\Sht^{(G,\mu)}_{(H_1,\l),(H_2,\k),K_{0}gK_{0}}$, should also be the intersection number of $\Hk^{\l,\k}_{\cM,K_{0}gK_{0}}$ with the graph of Frobenius for $\cM^{G}_{H_{1}, H_{2},K_{0}gK_{0}}$. In other words, we are changing the order of intersection and leaving the ``Shtuka-like'' construction to the very last step. It is often true that $\Hk^{\l,\k}_{\cM,K_{0}gK_{0}}$ has the same dimension as $\cM^{G}_{H_{1}, H_{2},K_{0}gK_{0}}$, and the fundamental class of $\Hk^{\l,\k}_{\cM,K_{0}gK_{0}}$ induces an endomorphism of the cohomology of $\cM^{G}_{H_{1}, H_{2},K_{0}gK_{0}}$ which we denote by $[\Hk^{\l,\k}_{\cM,K_{0}gK_{0}}]$. The Lefschetz trace formula then gives the following heuristic identity
\begin{equation}\label{I Lef trace}
I^{(G,\mu)}_{(H_1,\l),(H_2,\k)}(\one_{K_{0}gK_{0}})''=''\Tr([\Hk^{\l,\k}_{\cM,K_{0}gK_{0}}]\circ\Frob, \cohoc{*}{\cM^{G}_{H_{1}, H_{2},K_{0}gK_{0}}\ot\kbar, \Ql}).
\end{equation}
The equal sign above is in quotation marks for at least two reasons: both sides may diverge;  change of the order of intersection  needs to be justified.

Let $\cB^{G}_{H_{1}, H_{2},K_{0}gK_{0}}$ be the Hitchin base and $h^{G}_{H_{1},H_{2}}$  be the Hitchin map as in \eqref{rel Hitchin}. Observe that for various $0\le i\le r$, the compositions  $h^{G}_{H_{1}, H_{2}}\circ m_{i}: \Hk^{\l,\k}_{\cM,K_{0}gK_{0}}\to\cB^{G}_{H_{1}, H_{2},K_{0}gK_{0}}$ are all the same. On the other hand, the Frobenius of $\cM^{G}_{H_{1}, H_{2},K_{0}gK_{0}}$ covers the Frobenius of $\cB^{G}_{H_{1}, H_{2},K_{0}gK_{0}}$. Therefore diagram \eqref{alt ShtM} induces a map
\begin{equation*}
\Sht^{(G,\mu)}_{(H_1,\l),(H_2,\k),K_{0}gK_{0}}\to \cB^{G}_{H_{1}, H_{2},K_{0}gK_{0}}(k)
\end{equation*}
which simply says that $\Sht^{(G,\mu)}_{(H_1,\l),(H_2,\k),K_{0}gK_{0}}$ decomposes into a disjoint union
\begin{equation}\label{decomp ShtM}
\Sht^{(G,\mu)}_{(H_1,\l),(H_2,\k),K_{0}gK_{0}}=\coprod_{a\in \cB^{G}_{H_{1}, H_{2},K_{0}gK_{0}}(k)}\Sht^{(G,\mu)}_{(H_1,\l),(H_2,\k),K_{0}gK_{0}}(a).
\end{equation}
The action of $[\Hk^{\l,\k}_{\cM,K_{0}gK_{0}}]$ on the cohomology of $\cM^{G}_{H_{1}, H_{2},K_{0}gK_{0}}$ can be localized to an action on the complex $\bR h^{G}_{H_{1}, H_{2},!}\Ql$ using the formalism of cohomological correspondences. Then we may rewrite \eqref{I Lef trace} as
\begin{equation}\label{I Lef}
I^{(G,\mu)}_{(H_1,\l),(H_2,\k)}(\one_{K_{0}gK_{0}})''=''\sum_{a}\Tr([\Hk^{\l,\k}_{\cM,K_{0}gK_{0}}]_{a}\circ\Frob_{a}, (\bR h^{G}_{H_{1}, H_{2},!}\Ql)_{a}),
\end{equation}
where $a$ runs over $\cB^{G}_{H_{1}, H_{2},K_{0}gK_{0}}(k)$.

Comparing \eqref{I Lef} with \eqref{RTr Lef}, we see the only difference is the insertion of the operator  $[\Hk^{\l,\k}_{\cM,K_{0}gK_{0}}]_{a}$ acting on the stalk $(\bR h^{G}_{H_{1}, H_{2},!}\Ql)_{a}$.

\begin{exam}\label{ex:T cont}  Suppose $G=\PGL_{2}$ and $H_{1}=H_{2}=T$ as in Example \ref{ex:T}. Let $r\ge0$ be even. We pick $\l,\k\in\{\pm1\}^{r}$ with total sum zero, and form $\Sht^{\l}_{T}$ and $\Sht^{\k}_{T}$. Let $\mu=(\mu_{1},\cdots, \mu_{r})$ consist of minuscule coweights of $G$. In this situation we have $\dim \Sht^{\l}_{T}=\dim\Sht^{\k}_{T}=r=\frac{1}{2}\dim \Sht^{\mu}_{G}$. Below we give more explicit descriptions of $\cM^{G}_{T,T,K_{0}gK_{0}}$ and $\Hk^{\l,\k}_{\cM,K_{0}gK_{0}}$. For simplicity, we assume that the double cover $\nu: X'\to X$ is \'etale.

As in Example \ref{ex:A}, it is more convenient to work with the test function $h_{D}$ rather than $\one_{K_{0}gK_{0}}$, where $D$ is an effective divisor on $X$. We denote the corresponding version of $\cM^{G}_{T,T,K_{0}gK_{0}}$ by $\cM^{G}_{T,T,D}$. To describe $\cM^{G}_{T,T,D}$, we first consider the moduli stack $\wt\cM^{G}_{T,T,D}$ classifying $(\cL, \cL', \ph)$ where $\cL$ and $\cL'$ are line bundles over $X'$, and $\ph: \nu_{*}\cL\to \nu_{*}\cL'$ is an injective map of coherent sheaves such that $\det(\ph)$ has divisor $D$. We then have $\cM^{G}_{T,T,D}=\wt\cM^{G}_{T,T,D}/\Pic_{X}$, where $\Pic_{X}$ acts by pulling back to $X'$ and simultaneously tensoring with $\cL$ and $\cL'$. 

It is more convenient to work with another description of $\cM^{G}_{T,T,D}$. The map $\ph: \nu_{*}\cL\to \nu_{*}\cL'$ is equivalent to the data of two maps
\begin{eqnarray*}
\a: \cL\to \cL', \quad  \b: \s^{*}\cL\to \cL'
\end{eqnarray*}
where $\s$ is the nontrivial involution of $X'$ over $X$. The determinant $\det(\ph)=\Nm(\a)-\Nm(\b)$, as sections of $\Nm(\cL)^{-1}\ot\Nm(\cL')$. Let $\cM^{G,\ds}_{T,T,D}\subset \cM^{G}_{T,T,D}$ be the open subset where $\a$ and $\b$ are nonzero. By recording the divisors of $\a$ and $\b$, we may alternatively describe $\cM^{G,\ds}_{T,T,D}$ as the moduli of pairs $(D_{\a}, D_{\b})$ of effective divisors on $X'$ of degree $d=\deg D$, such that there exists  a rational function $f$ on $X$ (necessarily unique) satisfying $\div(f)=\nu(D_{\a})-\nu(D_{\b})$ and $\div(1-f)=D-\nu(D_{\b})$.

The Hecke correspondence $\Hk^{\l,\k}_{\cM,D}$ is the composition of $r$ correspondences each of which is either $\cH_{+}$ or $\cH_{-}$ depending on whether $\l_{i}=\k_{i}$ or not. Over the open subset $\cM^{G,\ds}_{T,T,D}$,  $\cH_{+}$  can be described as follows: it classifies triples of effective divisors $(D_{\a}, D_{\b}, D'_{\b})$ on $X'$ such that $(D_{\a}, D_{\b})\in \cM^{G,\ds}_{T,T,D}$, and $D'_{\b}$ is obtained by changing one point of $D_{\b}$ by its image under $\s$. The two maps $p_{+}, q_{+}: \cH_{+}\to \cM^{G,\ds}_{T,T,D}$ send $(D_{\a}, D_{\b}, D'_{\b})$ to $(D_{\a}, D_{\b})$ and $(D_{\a}, D'_{\b})$. Similarly, over $\cM^{G,\ds}_{T,T,D}$, $\cH_{-}$ classifies triples of effective divisors $(D_{\a}, D'_{\a}, D_{\b})$ on $X'$ such that $D'_{\a}$ is obtained by changing one point of $D_{\a}$ by its image under $\s$. 
\end{exam}

\subsection{Application to $L$-functions}\label{ss:YZ}
In the work \cite{YZ1}, we considered the case $G=\PGL_{2}$ and the moduli of Shtukas $\Sht^{r}_{G}$ without level structures, where $r$ stands for the $r$-tuple $\mu=(\mu_{1},\cdots, \mu_{r})$ consisting of minuscule coweights of $G$ (so $r$ is even).  Let $\nu:X'\to X$ be an unramified double cover. The Heegner-Drinfeld cycle we considered was the one introduced in Example \ref{ex:T}, i.e., $\Sht^{\l}_{T}$ for  $\l\in\{\pm1\}^{r}$. We consider the lifting of the natural map  $\th:\Sht^{\l}_{T}\to \Sht^{r}_{G}$
\begin{equation*}
\th': \Sht^{\l}_{T} \to \Sht'^{r}_{G}:=\Sht^{r}_{G}\times_{X^{r}}X'^{r}.
\end{equation*} 
Since $\Sht^{\l}_{T}$ is proper of dimension $r$, the Heegner-Drinfeld cycle $\cZ^{\l}_{T}:=\th'_{*}[\Sht^{\l}_{T}]$ is an $r$-dimensional proper cycle in the $2r$-dimensional $\Sht'^{r}_{G}$. Therefore $\cZ^{\l}_{T}$ defines a class  $Z^{\l}_{T}\in\cohoc{2r}{\Sht'^{r}_{G}\ot\kbar,\Qlbar}(r)$. 

Now let $\pi$ be an everywhere unramified cuspidal automorphic representation of $G(\AA)$ with coefficients in $\Qlbar$. By the coarse cohomological spectral decomposition for $\cohoc{2r}{\Sht'^{r}_{G}\ot\kbar,\Qlbar}$ (see Theorem \ref{th:coho spec}), we may project $Z^{\l}_{T}$ to the $\chi_{\pi}$-isotypical summand, and denote the resulting class by $Z^{\l}_{T,\pi}\in\cohoc{2r}{\Sht'^{r}_{G}\ot\kbar,\Qlbar}[\chi_{\pi}]$.

\begin{theorem}[\cite{YZ1}]\label{th:YZ1} We have
\begin{equation*}
\jiao{Z^{\l}_{T,\pi}, Z^{\l}_{T,\pi}}_{\Sht'^{r}_{G}}=\frac{q^{2-2g}}{2(\log q)^{r}}\frac{\sL^{(r)}(\pi_{F'}, 1/2)}{L(\pi,\Ad,1)}
\end{equation*}
where
\begin{itemize}
\item $\jiao{Z^{\l}_{T,\pi}, Z^{\l}_{T,\pi}}_{\Sht'^{r}_{G}}$ is the self-intersection number of the cycle class $Z^{\l}_{T,\pi}$.
\item $\pi_{F'}$ is the base change of $\pi$ to $F'=k(X')$. 
\item $\sL(\pi_{F'},s)=q^{4(g-1)(s-1/2)}L(\pi_{F'},s)$ is the normalized $L$-function of $\pi_{F'}$ such that $\sL(\pi_{F'},s)=\sL(\pi_{F'},1-s)$.
\end{itemize}
\end{theorem}

In \cite{YZ2},  we extended the above theorem to allow the automorphic representation $\pi$ to have square-free level structures (which means the local representations $\pi_{v}$ are either unramified or an unramified twist of the Steinberg representation),  and to allow ramifications for the double cover $\nu:X'\to X$. We consider the moduli stack  $\Sht^{r}_{G}(\Sig;\Si)$ where $\Sig$ is a finite set of places where we add Iwahori level structures to the $G$-Shtukas; $\Si\subset \Sig$ is a subset of places where we impose supersingular conditions. The admissibility condition forces $r$ to have the same parity as $\#\Si$.

\begin{theorem}[\cite{YZ2}]\label{th:YZ2} Let $\pi$ be a cuspidal automorphic representation of $G(\AA)$ with square-free level $\Sigma$. Assume the double cover $\nu:X'\to X$ is unramified over $\Sigma$. Let $\Si\subset\Sigma$ be the places that are inert in $F'$. Let $r\in\ZZ_{\ge0}$ be of the same parity as $\#\Si$. Then for any $r_{1}, r_{2}\in\ZZ_{\ge0}$ such that $r_{1}+r_{2}=r$,  there is an explicit linear combination $Z^{r_{1}, r_{2}}_{T}$ of the cycles $\{Z^{\l}_{T}; \l\in\{\pm1\}^{r}\}$ such that 
\begin{equation*}
\jiao{Z^{r_{1},r_{2}}_{T,\pi}, Z^{r_{1}, r_{2}}_{T,\pi}}_{\Sht'^{r}_{G}(\Sig;\Si)}=\frac{q^{2-2g+\r/2-N}}{2(-\log q)^{r}}\frac{\sL^{(r_{1})}(\pi, 1/2)\sL^{(r_{2})}(\pi\otimes\eta_{F'/F}, 1/2)}{L(\pi,\Ad,1)}.
\end{equation*}
where
\begin{itemize}
\item $N=\deg\Sig$, and $\rho$ is the degree of the ramification locus of $\nu$.
\item $\sL(\pi,s)=q^{(2g-2+N/2)(s-1/2)}L(\pi,s)$ is the normalized $L$-function of $\pi$ such that $\sL(\pi,s)=\sL(\pi,1-s)$.
\item $\eta_{F'/F}$ is the  character of $F^{\times}\bs \AA^{\times}$ corresponding to the quadratic extension $F'/F$.   
\item $\sL(\pi\ot\eta_{F'/F},s)=q^{(2g-2+\rho+ N/2)(s-1/2)}L(\pi\ot\eta_{F'/F},s)$ is the normalized $L$-function of $\pi\ot\eta_{F'/F}$ such  that $\sL(\pi\ot\eta_{F'/F},s)=\sL(\pi\ot\eta_{F'/F},1-s)$.
\end{itemize}
\end{theorem}

When $r=0$, the above theorem is a special case of the Waldspurger formula \cite{W}, and our proof in this case is very close to the one given by Jacquet \cite{JW}.  When $r=1$ and $\#\Si=1$, the above theorem is an analogue of the Gross-Zagier formula (see \cite{GZ}) which expresses the first derivative of the base-change $L$-function of a cuspidal Hecke eigenform in terms of the height of Heegner points on the modular curve.  However our proof is very different from the original proof of the Gross-Zagier formula in that we do not need to explicitly compute either side of the formula.

\sss{Relation with the B-SD conjecture} Theorem \ref{th:YZ2} is applicable to those $\pi$ coming from semistable elliptic curves $E$ over the function field $F$. The relation of our result and the Birch--Swinnerton-Dyer conjecture for $E$ can be roughly stated as follows. Take $r$ to be the vanishing order of $L(E_{F'}, s)=L(\pi_{F'}, s-1/2)$ at $s=1$. According to the expectation \eqref{coho Sht}, $Z^{\l}_{T,\pi}$ is an element in $\pi^{K}\ot \cohog{1}{X'\ot\kbar, j_{!*}\nu^{*}\rho_{\pi}}^{\otimes r}$. The $2$-dimensional $\ell$-adic Galois representation $\rho_{\pi}$ attached to $\pi$ is the Tate module of $E$, therefore $L(E_{F'},s)=\det(1-q^{-s}\Frob|\cohog{1}{X'\ot\kbar, j_{!*}\nu^{*}\rho_{\pi}})$. The standard conjecture predicts that the Frobenius acts semisimply on $\cohog{1}{X'\ot\kbar, j_{!*}\nu^{*}\rho_{\pi}}$, hence the multiplicity of the Frobenius eigenvalue $q$ should be $r$. We expect $Z^{\l}_{T,\pi}$ to lie in $\pi^{K}\ot \wedge^{r}(\cohog{1}{X'\ot\kbar, j_{!*}\nu^{*}\rho_{\pi}}^{\Fr=q})$, and giving a basis for this hypothetically $1$-dimensional space. However, currently we do not have a way to construct rational points on $E$ from the Heegner-Drinfeld cycle $\Sht^{\l}_{T}$.

\sss{} The method to prove Theorems \ref{th:YZ1} is by comparing the Shtuka version of the relative trace $I^{G}_{T,T}(f)=\jiao{\cZ^{\l}_{T}, f*\cZ^{\l}_{T}}_{\Sht'^{r}_{G}}$ as in \eqref{defn I} with the usual relative trace of the kind in Example \ref{ex:A}. More precisely, for the triple $(G,A,A)$ considered in Example \ref{ex:A}, we consider the relative trace involving a complex variable $s$
\begin{equation*}
\RTr^{G}_{A,(A,\eta)}(f,s)=\jiao{\ph_{!}|\cdot|^{s}_{A(F)\bs A(\AA)}, \ph_{!}(\eta|\cdot|^{s})_{A(F)\bs A(\AA)}}_{L^{2}(G(F)\bs G(\AA), \mu_{G})}.
\end{equation*}
Here $|\cdot|:A(F)\bs A(\AA)=F^{\times}\bs \AA^{\times}\to q^{\ZZ}$ is the global absolute value function, and $\y=\y_{F'/F}$. Let $J_{r}(f)$ be the $r$th derivative of $\RTr^{G}_{A,(A,\eta)}(f,s)$ at $s=0$.  The key to the proof is to establish the following identity of relative traces for all spherical Hecke functions $f$
\begin{equation}\label{I=J}
I^{G}_{T,T}(f)=(\log q)^{-r}J_{r}(f).
\end{equation}

To prove this identity, it suffices to consider $f=h_{D}$ for effective divisors $D$ on $X$  (see Example \ref{ex:A}). The moduli stacks $\cM^{G}_{A,A,D}$ in Example \ref{ex:A} and $\cM^{G}_{T,T,D}$ in Example \ref{ex:T cont} share the same Hitchin base $\cB_{D}=\cohog{0}{X,\cO_{X}(D)}$. We may fix a degree $d$ and let $D$ vary over effective divisors of degree $d$ and get Hitchin maps $f_{d}: \cM^{G}_{A,A,d}\to \cB_{d}$ and $g_{d}: \cM^{G}_{T,T,d}\to \cB_{d}$.  Formulae \eqref{RTr Lef} and \eqref{I Lef} suggest that we should try to prove an identity between the direct image complexes of $f_{d}$ and $g_{d}$. The new geometric input here is the action of the Hecke correspondences $[\Hk^{\l,\l}_{\cM_{d}}]$ on the complex $\bR g_{d!}\Ql$, which is the $r$-th iteration of the action of the correspondence $[\cH_{+}]$ defined in Example \ref{ex:T cont}. It turns out that the eigenvalues of the action of $[\cH_{+}]$ on $\bR g_{d!}\Ql$ match exactly with the factors coming from taking the derivative of the relative trace $\RTr^{G}_{A,(A,\eta)}(f,s)$, which explains why derivatives of automorphic quantities are indeed geometric.

\subsection{Arithmetic fundamental lemma}\label{ss:AFL}
Generalizing Theorems \ref{th:YZ1} and \ref{th:YZ2} to higher rank groups would involve intersecting non-proper cycles in an ambient stack which is not of finite type. This is the same issue as the non-convergence of the naive relative trace \eqref{rel trace}, therefore a certain truncation and regularization procedure is needed. There is, however, a local version of such results that can be proved for higher rank groups. One example of such a local version is the Arithmetic Fundamental Lemma formulated by W. Zhang \cite{ZAFL} originally for Rapoport-Zink spaces.  In \cite{YAFL}, we stated a higher derivative extension of W.Zhang's conjecture in the function field case, and sketched a proof.  This was the first time higher derivatives of automorphic quantities were related to geometry, and it partially motivated the later work \cite{YZ1}.

\sss{Local Shtukas} The moduli of Shtukas has a local version. Fix a local function field $F_{x}$ with ring of integers $\cO_{x}$. In the diagram  \eqref{defn ShtG} defining the moduli of Shtukas, we may replace $\Bun_{G}$ by the affine Grassmannian $\Gr_{G}$, and replace $\Hk^{\mu}_{G}$ by an iterated Hecke correspondence for $\Gr_{G}$ over a formal disk $\D_{r}=\Spf (\cO_{x}\hotimes \cdots \hotimes \cO_{x})$ of dimension $r$. We also have the freedom of changing the Frobenius morphism on $\Gr_{G}$ by the composition $b\circ \Fr: \Gr_{G}\to \Gr_{G}$, where $b$ is an element of the loop group of $G$ giving the datum of a $G$-isocrystal. The resulting object ${}^{b}\Sht^{\mu, \loc}_{G}$ by forming the Cartesian square as \eqref{defn ShtG} is called the moduli space of local Shtukas, and it is a formal scheme over $\D_{r}$. The special fiber of ${}^{b}\Sht^{\mu, \loc}_{G}$ is an iterated version of the affine Deligne-Lusztig variety. The twisted centralizer $G_{b}$ of $b$ is an inner form of a Levi subgroup of $G$, and $G_{b}(F_{x})$ acts on ${}^{b}\Sht^{\mu,\loc}_{G}$.

For a subgroup $H\subset G$ with a sequence of coweights $\l=(\l_{1},\cdots, \l_{r})$ such that $\l_{i}\le_{G}\mu_{i}$ for all $i$ and an $H$-isocrystal $b_{H}$ compatible with $b$, we have a morphism $\th^{\loc}: {}^{b_{H}}\Sht^{\l,\loc}_{H}\to {}^{b}\Sht^{\mu,\loc}_{G}$ over $\D_{r}$. This morphism is a closed embedding because $\Gr_{H}\to \Gr_{G}$ is. We define the local Heegner-Drinfeld cycle ${}^{b_{H}}\cZ^{\l,\loc}_{H}$ as the image of $\th^{\loc}$.

If we have two local Heegner-Drinfeld cycles ${}^{b_{1}}\cZ^{\l,\loc}_{H_{1}}$ and ${}^{b_{2}}\cZ^{\k, \loc}_{H_{2}}$ in ${}^{b}\Sht^{\mu,\loc}_{G}$ with  complementary dimensions, and if $\mu_{i}$ are minuscule and the reduced structure of their intersection is proper over $k$, we may ask for their intersection number in ${}^{b}\Sht^{\mu,\loc}_{G}$. More generally, if $\d\in G_{b}(F_{x})$, we may consider the intersection number
\begin{equation*}
I_{\d}=\jiao{{}^{b_{1}}\cZ^{\l,\loc}_{H_{1}}, \d\cdot {}^{b_{2}}\cZ^{\k,\loc}_{H_{2}}}_{{}^{b}\Sht^{\mu,\loc}_{G}}
\end{equation*}
using the action of $G_{b}(F_{x})$ on ${}^{b}\Sht^{\mu,\loc}_{G}$. When $\mu=0$, this is the same as the local orbital integral $J^{G}_{H_{1},H_{2},\d}(\one_{G(\cO_{x})})$ (see \eqref{rel orb int}) for the relative trace formula of the triple $(G,H_{1}, H_{2})$.

\begin{exam} Let $F'_{x}/F_{x}$ be an unramified quadratic extension, with ring of integers  $\cO'_{x}$ and residue field $k'_{x}$. Fix a Hermitian vector space $W_{n,x}$ of dimension $n$ over $x$, and let $\Ug_{n}$ be the unitary group of $W_{n,x}$. We define the moduli of local Shtukas $\Sht^{r, \loc}_{\Ug_{n}}$ over $\D'_{r}=\Spf(\cO'_{x}\hotimes_{k'_{x}}\cdots\hotimes_{k'_{x}}\cO'_{x})$ in the following way. Let $\Gr_{\Ug_{n}}$ be the  affine Grassmannian classifying self-dual lattices in $W_{n,x}$. Since $\Ug_{n}$ is split over $F'_{x}$, the base change $\Gr_{\Ug_{n}}\ot_{k_{x}}k'_{x}$ can be identified with the affine Grassmannian $\Gr_{\GL_{n}}\ot_{k_{x}}k'_{x}$ classifying $\cO'_{x}$-lattices in $W_{n,x}$.  We have the local Hecke correspondence $\Hk^{\loc}_{\Ug_{n}}$ over $\D'_{1}$ which, after identifying $\Gr_{\Ug_{n}}\ot_{k_{x}}k'_{x}$ with $\Gr_{\GL_{n}}\otimes_{k_{x}}k'_{x}$, corresponds to the upper modification of lattices in $W_{n,x}$ of colength one. Let $\Hk^{r,\loc}_{\Ug_{n}}$ be the $r$-fold composition of $\Hk^{\loc}_{\Ug_{n}}$ as a correspondence, so $\Hk^{r,\loc}_{\Ug_{n}}\to \D'_{r}$. Then $\Sht^{r, \loc}_{\Ug_{n}}$ is defined using the Cartesian diagram
\begin{equation*}
\xymatrix{      \Sht^{r, \loc}_{\Ug_{n}}\ar[r]\ar[d] &  \Hk^{r,\loc}_{\Ug_{n}}  \ar[d]  \\
\Gr_{\Ug_{n}}\ar[r]^-{(\id, \Fr)} & \Gr_{\Ug_{n}}\times \Gr_{\Ug_{n}}}
\end{equation*}
The group $\Ug_{n}(F_{x})$ acts on $\Sht^{r, \loc}_{\Ug_{n}}$. We remark that  $ \Sht^{r, \loc}_{\Ug_{n}}$ is formally smooth over $\D'_{r}$ of relative dimension $r(n-1)$.

Now suppose we are in the situation of Example \ref{ex:JR Ug} so we have an orthogonal decomposition $W_{n,x}=W_{n-1,x}\oplus F'_{x}e_{n}$ and $(e_{n},e_{n})=1$.  Adding a standard lattice $\cO'_{x}e_{n}$ gives an embedding $\Gr_{\Ug_{n-1}}\incl \Gr_{\Ug_{n}}$ compatible with the Hecke modifications, hence induces an embedding $\Sht^{r, \loc}_{\Ug_{n-1}}\incl \Sht^{r, \loc}_{\Ug_{n}}$. Consider the diagonal map
\begin{equation*}
\D: \Sht^{r, \loc}_{\Ug_{n-1}}\to \Sht^{r, \loc}_{\Ug_{n-1}}\times_{\D'_{r}}\Sht^{r, \loc}_{\Ug_{n}}.
\end{equation*}
The image of $\D$ gives an $r(n-1)$-dimensional cycle $\cZ^{r,\loc}_{x}$ in the $2r(n-1)$-dimensional ambient space $\Sht^{r, \loc}_{\Ug_{n-1}}\times_{\D'_{r}}\Sht^{r, \loc}_{\Ug_{n}}$.  For  strongly regular semisimple $\d\in \Ug_{n}(F_{x})$ (with respect to the conjugation action by $\Ug_{n-1}$)  we form the intersection number
\begin{equation*}
I^{\Ug,r}_{x,\d}=\jiao{\cZ^{r,\loc}_{x}, (\id\times\d)\cZ^{r,\loc}_{x}}_{\Sht^{r,\loc}_{\Ug_{n-1}}\times_{\D'_{r}}\Sht^{r,\loc}_{\Ug_{n}}}.
\end{equation*}
This makes sense because the support of the intersection of the two cycles $\cZ^{r,\loc}_{x}$ and $(\id\times\d)\cZ^{r,\loc}_{x}$ is proper. When $r=0$, we have $I^{\Ug,r}_{x,\d}$ is either equal to $J^{\Ug}_{x,\d}(\one_{\Ug(\L_{n,x})})$ as in Example \ref{ex:JR Ug} if a self-dual lattice $\L_{n,x}\subset W_{n,x}$ exists, or $0$ otherwise.
\end{exam}

To state the higher arithmetic fundamental lemma, we introduce a variant of the orbital integral \eqref{orb int GL} with a complex variable $s$
\begin{equation*}
J^{\GL}_{x,\g}(f,s)=\int_{\GL_{n-1}(F_{x})}f(h^{-1}\g h)\eta_{x}(\det h)|\det h|^{s}dh,\quad \g\in \Sg_{n}(F_{x}), f\in C_{c}^{\infty}(\Sg_{n}(F_{x})).
\end{equation*}
For $\g$ strongly regular semisimple, $J^{\GL}_{x,\g}(f,s)$ is a Laurent polynomial in $q^{s}_{x}$, where $q_{x}=\#k_{x}$. Let
\begin{equation*}
J^{\GL,r}_{x,\g}(f)=\left(\frac{d}{ds}\right)^{r}|_{s=0}J^{\GL}_{x,\g}(f,s).
\end{equation*}

\begin{theorem}[\cite{YAFL}] Let $\g\in \GL_{n}(F_{x})$ and $\d\in \Ug_{n}(F_{x})$ be strongly regular semisimple with the same invariants in the sense of  Examples \ref{ex:JR GL} and \ref{ex:JR Ug}. Then
\begin{equation}\label{rough AFL}
I^{\Ug,r}_{x,\d}=c(\log q_{x})^{-r}J^{\GL,r}_{x,\g}(\one_{\Sg_{n}(\cO_{x})})
\end{equation}
with an explicit constant $c$ depending on $r$ and the invariants of $\g$.
\end{theorem}
The proof consists of the following main steps.
\begin{enumerate}
\item Prove a global analogue of \eqref{rough AFL}. Consider the triple $(G,H_{1}, H_{2})$ as in Example \ref{ex:JR Ug} with coweights $\mu$ for $G$ and $\l$ for $H_{1}=H_{2}$ being the first fundamental coweights. Let the triple $(G',H'_{1}, H'_{2})$ be as in Example \ref{ex:JR GL}. Recall that for the global situation, the intersection number $I^{(G,\mu)}_{(H_{1},\l),(H_{2},\l)}(\one_{K_{0}gK_{0}})$ is the degree of a certain $0$-cycle on $\Sht^{(G,\mu)}_{(H_{1},\l),(H_{2},\l),K_{0}gK_{0}}$ introduced in the diagram \eqref{defn ShtM}. On the other hand, we have a decomposition \eqref{decomp ShtM} of $\Sht^{(G,\mu)}_{(H_{1},\l),(H_{2},\l),K_{0}gK_{0}}$ into a disjoint union of $\Sht^{(G,\mu)}_{(H_{1},\l),(H_{2},\l),K_{0}gK_{0}}(a)$ indexed by $k$-points of the base $\cB^{G}_{H_{1},H_{2},K_{0}gK_{0}}$. For strongly regular semisimple $a$, $\Sht^{(G,\mu)}_{(H_{1},\l),(H_{2},\l),K_{0}gK_{0}}(a)$ is proper for any $g$, so we can talk about the degree of the $a$-component of the zero cycle $\cZ^{\l}_{H_{1}}\cdot (f*\cZ^{\l}_{H_{2}})$, denoted $\jiao{\cZ^{\l}_{H_{1}}, f*\cZ^{\l}_{H_{2}}}_{a}$. The global analogue of \eqref{rough AFL} means proving an identity of the form \eqref{I=J}, but with both sides replaced by their $a$-components. One can prove such a global identity by analyzing the direct image complexes of the Hitchin maps $h^{G}_{H_{1}, H_{2}}$ and $h^{G'}_{H'_{1}, H'_{2}}$ using sheaf-theoretic methods, as we did in \cite{YJR} and \cite{YZ1}.

\item Deduce the arithmetic fundamental lemma from the global identity. The moduli $\Sht^{\mu,\loc}_{G}$ of local Shtukas for $G$ is related to a formal completion of the global moduli stack $\Sht^{\mu}_{G}$ by a uniformization diagram, analogous to the one relating Rapoport-Zink spaces and Shimura varieties.  Using the uniformization, one can express $\jiao{\cZ^{\l}_{H_{1}}, f*\cZ^{\l}_{H_{2}}}_{a}$ as a finite sum, where each summand is a product of usual orbital integrals and intersection numbers of the form $I^{\Ug, r_{i}}_{x_{i},\d}$ (with $\sum_{}r_{i}=r$). There is a similar product formula for the global orbital integral for $(G',H'_{1}, H'_{2})$. By choosing $a$ appropriately we may deduce the local identity \eqref{rough AFL} from the global one using the product expansions and the known fundamental lemma (Theorem \ref{th:JR}).

\end{enumerate}

\end{document}